\newcommand{\bel}{\begin{equation}\label}
\newcommand{\ee}{\end{equation}}

\documentclass[11 pt,letterpaper]{amsart}
\usepackage{graphicx,epsfig}
\usepackage{amssymb}
\usepackage{todonotes}

\newenvironment{namelist}[1]{%
\begin{list}{}
     {
      
      \settowidth{\labelwidth}{#1}
      \setlength{\leftmargin}{1.1\labelwidth}
               }
      }{%
\end{list}}
\newenvironment{Example}{\begin{Ex}\rm}{\smallskip\end{Ex}}
\newtheorem{thm}{Theorem}

\newtheorem{Ex}{Example}[section]
\newtheorem{defn}{Definition}
\newtheorem{rem}{Remark}
\begin{document}
\title[Extremal Kendall convolution]{How exceptional is the extremal Kendall and Kendall-type convolution}
\author[B.H. Jasiulis-Go{\l}dyn, J.K. Misiewicz, E. Omey, J. Weso{\l}owski]{Barbara H. Jasiulis-Go{\l}dyn$^1$, Jolanta K. Misiewicz $^2$, Edward Omey$^3$ and Jacek Weso{\l}owski$^4$}
\thanks{$^1$ Institute of Mathematics, University of Wroc{\l}aw, pl. Grunwaldzki 2, 50-384 Wroc{\l}aw, Poland, e-mail: barbara.jasiulis@math.uni.wroc.pl \\
$^{2, 4}$ Faculty of Mathematics and Information Science, Warsaw University of Technology, ul. Koszykowa 75, 00-662 Warsaw, Poland, \\
$^2$ e-mail: j.misiewicz@mini.pw.edu.pl, $^4$e-mail: j.wesolowski@mini.pw.edu.pl \\
$^3$ Faculty of Economics and Business-Campus Brussels, KU Leuven, Warmoesberg 26, 1000 Brussels, Belgium, e-mail:
edward.omey@kuleuven.be  
}
\begin{abstract} 
This paper deals with the generalized convolutions connected with the Williamson transform and the maximum operation. We focus on such convolutions which can  define transition probabilities of renewal processes. They should be monotonic since the described time or destruction does not go back, it should admit existence of a distribution with a lack of memory property because the analog of the Poisson process shall exist. Another valuable property is the simplicity of calculating and inverting the corresponding generalized characteristic function (in particular Williamson transform) so that the technique of generalized characteristic function can be used in description of our processes. The convex linear combination property (the generalized convolution of two point measures is the convex combination of several fixed measures), or representability (which means that the generalized convolution can be easily written in the language of independent random variables) - they also facilitate the modeling of real processes in that language.

We describe examples of generalized convolutions having the required properties ranging from the maximum convolution and its simplest generalization - the Kendall convolution (associated with the Williamson transform), up to the most complicated here - Kingman convolution. It is novel approach to apply in the extreme value theory. Stochastic representation of the Kucharczak-Urbanik in the order statistics terms is proved, which open new paths to investigate Archimedean copulas.

This paper open the door to solve an old open problem of the relationship between copulas and generalized convolutions mentioned by B. Schweizer and A. Sklar in 1983. This indicates the path of further research towards extremes and dependency modelling.
\end{abstract}
\maketitle

\tableofcontents
\section{Motivations}
The most important motivation to write this paper is to solve the open problem mentioned in \cite{Schweizer_Sklar} (see problem 7.9.9, p.122) to make a measurable calibration of copulas in the context of generalized convolutions, in particular extremes. 
We present here the first step towards solving this problem and present the characterization of the maximum convolution, responsible for the construction of extremes, in the terms of other convolutions, which in their algebras already have well-developed tools (e.g. the Williamson transform corresponding to the Kendall convolution).
Moreover, we are building new tools for the development of the L\'evy processes theory in generalized sense (\cite{BJMR}), in particular  extremes, and a groundbreaking renewal theory such that renewal functions have clear analytical formulas (see \cite{renewal}).
\vspace{5mm}

Real world phenomena are often effects of accumulation processes. The most natural description of a accumulation process is through summation of components. However, accumulation of  
components can be described sometimes more adequately by the maximum 
function or by the $\ell_p$-norm of the vector of components, but often this dependence is more complicated and we are using some approximating methods. The origin of this point of view one can find in two papers: \cite{BJMR, renewal}, where the authors consider renewal theory and L\'evy processes in the generalized convolution sense. Following this papers, we propose here to use  generalized convolutions, in particular the Kendall convolution and the Kendall type convolutions to find connections with the extreme value theory. Our choice is motivated by many interesting properties of such convolutions including the close connection with the max-convolution and the extreme value theory, simplicity in calculating the corresponding characteristic functions and inverting these characteristic functions, representing convolutions by the convex linear combination of some measures or representing them by simple operations on some independent random variables. 
Since the Kendall convolution extends the concept of the max-convolution we call it extremal Kendall convolution to emphasize this property.

In this paper one can find the precise description of Kendall convolution and the Kendall-type convolutions, their exceptional properties and applications in stochastic models - some of them have not been known yet. We analyze also some examples of convolutions with similar properties. Finally, we prove a stochastic representation of the Kucharczak-Urbanik convolution in the terms of order statistics. This is the starting point for a novel approach to Archimedean copulas which are now extensively studied (\cite{CCEHLJ,Larsson, McNeilNesholva1, McNeilNesholva2}).\\

Generalized convolutions were defined and intensively studied by K. Urbanik in \cite{Urbanik64, Urbanik73, Urbanik76, Urbanik84, Urbanik86}. His work has its origin in the paper of Kingman \cite{King}, where the first generalized convolution, called now the Kingman (or Bessel) convolution, was defined. This convolution - the ancestor of all generalized convolutions - is strictly connected with a Wiener process in $\mathbb{R}^n$ and the Bessel process describing the distance of the walking particle from the origin. 

For a while it was not clear that the class of generalized convolutions is rich enough to be interesting for applications and useful in stochastic simulation and mathematical modeling, but by now we know that this class is very rich, worth studying. It turned out, for example, that each generalized convolution has its own Gaussian distribution, exponential law and Poisson process with corresponding distribution with lack of memory property (see \cite{Poisson, factor, weaklevy, Kaniowski}). The origin of some generalized convolutions one can find also in Delphic semi-groups (\cite{Gilewski1, Kendall}). A different approach to generalized convolutions appeared in the theory of harmonic analysis, see e.g. \cite{Levitan, Sousa1(19), Sousa2(19)}. \\

The classical convolution, corresponding to the summation of independent random variables and the max-convolution corresponding to taking the maximum of independent random variables,  are examples of generalized convolutions. The extreme value theory described e.g. in \cite{Embrechts,Salvadori}, based on the max-convolution, is widely developed and is applied e.g. in modelling rare events with important consequences, like floods, hurricanes  (see \cite{Asmussen, Chen, Salvadori, Toulemonde}).  We focus here on the Kendall convolution, defined by Urbanik in \cite{Urbanik88} which  can be used to model e.g. some hydrological phenomena: pretty stable behaviour of the "natural" water level together with rarely appearing floods. We describe some  distributional properties of the Kendall and Kendall type convolutions (see \cite{KendallWalk,KucharczakUrbanik, Mej18l-c}).  Especially interesting and useful in modeling extremal events is that for Kendall and Kendall type convolutions the convolution of two measures with compact supports can have heavy tail. \\

In Section 2 we present  basics of the theory of generalized convolutions.  \\

Section 3 contains a list of generalized convolutions studied in this paper. \\

For Section 4 let us remind first that each generalized convolution corresponds to its own integral transform, for details and basic properties see \cite{Bin71, Bin84, Urbanik73, Urbanik84, Urbanik86, Urbanik88, Urbanik93}. We describe some properties of the Kendall convolution through  its generalized characteristic function - the Williamson transform. Especially simple and clear here is the inversion formula. More information and details one can find in \cite{BJMR} or \cite{Williamson}. The Williamson transform is also used in copula theory (see e.g. \cite{McNeilNesholva1, McNeilNesholva2}) since it is a generator of Archimedean copulas. For asymptotic properties of the Williamson transform see \cite{MBE}, \cite{renewal} and \cite{Larsson}. In Subsection 4.1 we draw the reader attention to the fact that the generalized convolution can be defined by the corresponding integral transform as the proper generalized characteristic function. It turned out that such approach was already considered in the area of Harmonic Analysis and Theory of Special Functions, see e.g. \cite{Levitan, Sousa1(19), Sousa2(19)}.
However generalized convolutions considered in those papers may not satisfy all Urbanik's assumptions. \\

In Section 5  we show that for $\alpha \leqslant 1$ there exists a (weakly stable) distribution $\mu$  such that the Kendall convolution $\lambda_1 \vartriangle_{\alpha} \lambda_2$  can be defined by the following equation:
\begin{equation}\label{astt}
\bigl( \lambda_1 \vartriangle_{\alpha} \lambda_2 \bigr) \circ \mu = \bigl( \lambda_1 \circ \mu\bigr) \ast \bigl( \lambda_2 \circ \mu\bigr),
\end{equation}
where $\ast$ is the classical convolution and the operation $\circ \colon \mathcal{P}_+^2 \rightarrow \mathcal{P}_+$  is defined as follows: $\mathcal{L}(\theta_1) \circ \mathcal{L}(\theta_2) = \mathcal{L}(\theta_1 \, \theta_2)$ for independent random variables $\theta_1, \theta_2$. Generalized convolutions with this property are called weak generalized convolutions. We indicate which of the convolutions which we considered  are weak.  \\

In Section 6 we study properties of generalized convolution allowing the construction of the corresponding Poisson process. We start from the monotonicity property stating that the generalized sum of positive random variables cannot be smaller than their maximum - this is necessary to have positive increments (of time). Not every generalized convolution has this property. We also study existence of distributions with the lack of memory property with respect to a given generalized convolution. 
For some convolutions such distributions do not exists. The main result of this section, Theorem 5, gives a few equivalent conditions for monotonic convolution to allow the existence of a distribution with the lack of memory property. We indicate such convolutions among ones we consider in this paper, e.g. for the Kendall convolution,
$\triangle_{\alpha}$, the power distribution ${\rm pow}(\alpha)$ with the density $\alpha x^{\alpha -1} \mathbf{1}_{[0,1]}(x)$ has the lack of memory property. \\

In Section 7 we show that for the Kendall
convolution, $\vartriangle_{\alpha}$, for $\alpha\le 1$ there exists a distribution $\nu$, which is weakly stable with respect to max-convolution, such that for any $\lambda_1,\lambda_2\in \mathcal P_+$ (probability
measures on $[0,\infty))$ we have
$$
\bigl( \lambda_1 \vartriangle_{\alpha} \lambda_2 \bigr) \circ \nu = \bigl( \lambda_1 \circ \nu\bigr) \triangledown \bigl( \lambda_2 \circ \nu\bigr),
$$
with the max-convolution $\triangledown$ defined by $\mathcal{L}(\theta_1) \triangledown \mathcal{L}(\theta_2) = \mathcal{L} ( \max\{ \theta_1, \theta_2\})$, where $\theta_1$ and $\theta_2$ are independent positive random variables. We have also the following property which, as it will be shown in Section 3, trivially follows from the definition of the Kendall convolution:  $1 = \left(\delta_a \vartriangle_{\alpha} \delta_b\right) ([\max\{a,b \}, \infty))$ $ > \left(\delta_a  \vartriangle_{\alpha} \delta_b\right) ((\max\{a,b \}, \infty))$.   By these properties we can model such processes as the change of water level in the river in the continuous time which is  pretty stable most of the time but sometimes goes into extremes. \\

An equivalent definition of Kendall convolution presented in Section 8, states that the Kendall convolution of two Dirac measures, $\delta_a$, $\delta_b$, is a convex linear combination of two fixed measures with coefficients of this combination depending on $a$ and $b$. In \cite{MisJas2} it was shown that the Kendall convolution is the only generalized convolution with this property. It was shown in \cite{Mej18l-c} that if the generalized convolution of $\delta_a$ and $\delta_b$ is a convex combination of $n$ fixed measures and with coefficients of this combination depending on $a$ and $b$ then the generalized convolution is similar to the Kendall convolution. We call them the Kendall-type convolutions. Such convex combination properties are not only useful in explicit calculations, but they allow to define  a family of  integral transforms parametrized by $n\ge 2$ extending in this way the Williamson transform (which covers the case $n=2$). \\

Finally, in Section 9 we focus on preparation for studying path properties of the L\'evy processes with respect to generalized convolution. In order to make it possible we need to express the given convolution in the language of operations on independent random variables. Such a construction for a given generalized convolution  is called representability (for details see \cite{BJMR}). Here we study a simplified version of this property expressing a generalized convolution of two measures $\lambda_1 \diamond \lambda_2$ corresponding to the independent random variables $\theta_1, \theta_2$ as a distribution of an explicitly defined variable $\Psi(\theta_1, \theta_2)$. If $\Psi(\theta_1, \theta_2)(\omega) = \overline{\Psi}(\theta_1(\omega), \theta_2(\omega))$ a.e. for some measurable function $\overline{\Psi}\colon \mathbb{R}^2 \rightarrow \mathbb{R}$, then $\overline{\Psi}(x,y) = (x^p + y^p )^{1/p}$ for some $p \in (0,\infty]$. In all other cases $\Psi(\theta_1, \theta_2)$ depends also on some other random variables. For example for the Kendall convolution we have $\mathcal{L}(\theta_1) \vartriangle_{\alpha} \mathcal{L}(\theta_2)$ is the distribution of 
$$
M \bigl(\mathbf{1}_{(\varrho^{\alpha},1]}(U) + \Pi_{2 \alpha} \mathbf{1}_{[0,\varrho^{\alpha}]} (U) \bigr),
$$
where $M = \max\{ \theta_1, \theta_2\}$, $\varrho = {{\min \{ \theta_1, \theta_2\}}/{\max \{ \theta_1, \theta_2\}}}$, $\Pi_q$ is a variable with the Pareto distribution $\pi_q$ and density $q x^{-q-1} \mathbf{1}_{[1,\infty)}(x)$,  $U$ has uniform distribution on $[0,1]$ and $\theta_1, \Theta_2, \Pi_{2\alpha}, U$ are independent.

{\bf Notation.}

Through this paper, by $\mathcal{P}_+$ (respectively $\mathcal{P}$ or more general $\mathcal{P}(\mathbb{E})$) we denote family of all probability measures on the Borel subsets of $\mathcal{\mathbb{R}_+} := [0, \infty)$ (respectively $\mathbb{R}$ or more general separable Banach space $\mathbb{E}$). The distribution of the random element $X$ is denoted by $\mathcal{L}(X)$.
A dilation family (rescalings) of operators $T_a:\mathcal{P}_+\to \mathcal{P}_+$, $a\in\mathbb{R}_+:=[0,\infty)$ is defined  for $\mu \in \mathcal{P}_+$ and any Borel set $B$ in the following way: $T_a\mu(B)=\mu(B/a)$ if $a>0$ and $T_0 \mu = \delta_0$. Equivalently, $T_a\mu = \mathcal{L}(aX)$ for $a \in \mathbb{R}_+$ and $\mathcal{L}(X) = \mu$.

\section{A primer on generalized convolutions}
The Kendall convolution is a well known  example of a generalized convolution defined by K. Urbanik in \cite{Urbanik64} and studied in \cite{Urbanik73, Urbanik76, Urbanik84, Urbanik86}. Urbanik was mainly interested in generalized convolutions on $\mathcal{P}_+$ and we shall do the same in this paper, but a wider approach is also possible. \\

In this section we present this part of the theory of generalized convolutions,  which is necessary for studying  properties of Kendall and other convolutions.  

\begin{defn}\label{geen}
A generalized convolution is a binary, associative and commutative operation $\diamond$ on $ \mathcal{P}_+$ with the following properties:
\begin{itemize}
\item[(i)] $\lambda \diamond \delta_0 = \lambda$ for all $\lambda \in \mathcal{P}_+$; \\
\item[(ii)] $(p\lambda_1 +(1-p)\lambda_2) \diamond \lambda = p(\lambda_1 \diamond \lambda) + (1-p)(\lambda_2 \diamond \lambda)$ for all $p \in [0,1]$ and $\lambda, \lambda_1, \lambda_2 \in \mathcal{P}_+$;\\
\item[(iii)] $T_a(\lambda_1 \diamond \lambda_2) = (T_a\lambda_1) \diamond (T_a\lambda_2)$ for all $a \geq 0$ and $\lambda_1, \lambda_2 \in \mathcal{P}_+$; \\
\item[(iv)] if $\lambda_n \Rightarrow \lambda$ and $\nu_n \Rightarrow \nu$, then $\lambda_n \diamond \nu_n \Rightarrow \lambda \diamond \nu$  for $\lambda_n, \mu_n, \lambda, \mu \in \mathcal{P}_+$, $n \in \mathbb{N}$, where $\Rightarrow $ denotes weak convergence;\\
\item[(v)] there exists a sequence of positive numbers $(c_n)$ and a probability measure $\nu\in\mathcal{P}_+$, $\nu\neq \delta_0$,  such that $T_{c_n} \delta_1^{\diamond n}\Rightarrow \nu$,  (here $\lambda^{\diamond n} = \underbrace{\lambda \diamond \lambda \diamond ... \diamond \lambda}_{n \;\mathrm{times}}$).
\end{itemize}
\end{defn}

\begin{rem}\label{x1def} Note that any generalized convolution $\diamond$ is uniquely determined by $\delta_x \diamond \delta_1$, $x\in[0,1]$. 
Indeed, by Def. \ref{geen},
\begin{itemize} \item first, for each choice of $a,b\in\mathbb{R}_+$ the measure $\delta_a \diamond \delta_b$ is uniquely determined by 
$$
\delta_a \diamond \delta_b = \left\{\begin{array}{ll} T_M \bigl( \delta_x \diamond \delta_1 \bigr), & \mathrm{if}\;M>0, \\
\delta_0, & \mathrm{if}\;M=0,\end{array}\right.
$$
where $M=a\vee b:=\max \{a,b\}$, $m =a\wedge b:= \min\{a,b\}$ and $x = \frac{m}{M}$; 
\item second, for  arbitrary measures  $\lambda_1, \lambda_2 \in \mathcal{P}_+$
$$
\lambda_1 \diamond \lambda_2= \int_0^{\infty} \int_0^{\infty} \left(\delta_a \diamond \delta_b \right) \, \lambda_1(da) \, \lambda_2 (db).
$$
\end{itemize}
\end{rem}
 
Characteristic functions are important tools for the analysis of classical convolution. It turns out that not  every generalized convolution allows a reasonable analog of characteristic function. The next definitions, introduced by K. Urbanik in \cite{Urbanik64}, select these convolutions for which such analog can be defined.

\begin{defn}
The class $\mathcal{P}_+$ equipped with the generalized convolution $\diamond$ is called a generalized convolution algebra and denoted by $(\mathcal{P}_+, \diamond)$. A continuous (in the sense of weak convergence of measures) mapping $h= h^{\diamond} \colon \mathcal{P}_+ \rightarrow \mathbb{R}$ is called a  homomorphism of the algebra  $(\mathcal{P}_+, \diamond)$ if for all $\lambda_1, \lambda_2 \in \mathcal{P}_+$  
\begin{itemize}
\item $h(\lambda_1 \diamond \lambda_2) = h(\lambda_1) h(\lambda_2)$ \\
and
\item $h(p\lambda_1 + (1-p)\lambda_2) = p h(\lambda_1) + (1-p) h(\lambda_2)$ for all $p \in [0,1]$.
\end{itemize}
Algebras admitting a non-trivial homomorphism  (i.e. $h\not\equiv 1$, $h \not\equiv 0$) and the corresponding generalized convolutions are called regular.
\end{defn}

\begin{defn}
For a regular algebra $(\mathcal{P}_+, \diamond)$ (or for the regular generalized convolution $\diamond$) we define a probability kernel $\Omega \colon \mathbb{R}_+ \rightarrow \mathbb{R}_+$ by
\begin{equation*}
\Omega(t) \overset{def}{=} h(T_t \delta_1) = h(\delta_t), \quad t \geqslant 0,
\end{equation*}
and a $\diamond$-generalized characteristic function $\Phi_{\lambda}^{\diamond}\colon  \mathbb{R}_+
\rightarrow \mathbb{R}_+$ of $\lambda\in\mathcal{P}_+$ as an integral transform with the kernel $\Omega$:
\begin{equation}\label{fidia}
\Phi_{\lambda}^{\diamond} (t) \overset{def}{=} \int_0^{\infty} \Omega(st) \lambda(ds) = h(T_t
\lambda),\quad t\in\mathbb{R}_+.
\end{equation}
\end{defn}

Note that if $X$ is a random variable with distribution $\lambda\in\mathcal{P}_+$ then 
$$
\Phi_{\lambda}^{\diamond} (t)={\bf E}\,\Omega(tX),\quad t\in\mathbb{R}_+.
$$
The function  $\Phi_{\lambda}^{\diamond}$ plays a similar role as the Laplace or Fourier transform for classical convolution on $\mathcal{P}_+$ or $\mathcal{P}$, respectively. Basic properties of $\diamond$-generalized characteristic functions are in \cite{weaklevy,Urbanik64}. For the present  paper it is important to know that each regular generalized convolution determines its generalized characteristic function uniquely up to a scale  constant. Moreover,  convergence of $\diamond$-generalized characteristic functions uniformly on  compact sets is equivalent to  weak convergence of the corresponding probability measures. \\

Some generalized convolutions admit only the existence of a function $h\colon \mathcal{P}_+ \rightarrow \mathbb{R}$  which has all the required properties of homomorphism except continuity. Equivalently, the corresponding probability kernel $\Omega \colon \mathbb{R}_+ \rightarrow \mathbb{R}$ is not continuous (and the corresponding generalized convolution is not regular).  For example max-convolution is not regular since it admits only one (up to a scale) probability kernel: $\Omega(x) = \mathbf{1}_{[0,1)}(x)$, which is evidently not continuous.  For such convolutions the corresponding generalized characteristic functions can be defined by \eqref{fidia}, but then some of  the properties, which hold in the regular case, may not be satisfied.

\section{Basic examples of generalized convolutions}
We present here a basic list of generalized convolutions defined uni\-que\-ly, according to Remark 1, by its value on $\delta_x \diamond \delta_y$ or $\delta_x \diamond \delta_1$ for $x \in (0,1)$. In the latter case, the values for $x \in \{ 0,1\}$ we define by continuity.  \\

{\bf Example 3.0.} The Kingman or Bessel convolution with parameter $s > -\frac{1}{2}$ is defined for $x,y \geqslant 0$ by 
$$
\delta_x \otimes_{\omega_s} \delta_y = \mathcal{L} \left( \sqrt{x^2 + y^2 + 2xy \theta_s} \right),
$$
where $\theta_s$ is a random variable with the following density function:
$$
f_s(t) = \frac{\Gamma(s+1)}{\sqrt{\pi} \Gamma( s + {1/2})} \, \bigl(1-t^2 \bigr)_+^{s-{1/2}},
$$
where $a_+ = \max\{a,0\}$.
The measure $\delta_x \otimes_{\omega_s} \delta_y$ has support $[ |x-y|, x+y ]$. If $n:= 2(s+1) \in \mathbb{N}$, $n >1$, then the variable $\theta_s$  can be identified  as $\theta_s = U_1$, where $\mathbf{U}_n = (U_1, \dots, U_n)$ is a random vector  having uniform distribution on the unit sphere $S_{n-1} \subset \mathbb{R}^n$. 

\begin{Example}
The classical convolution $\ast$ on $\mathcal{P}_+$ is given by
$$
\delta_x \ast \delta_y = \delta_{x+y}, \quad x,y \geqslant 0.
$$
\end{Example}

\begin{Example}
Symmetric convolution $\bowtie$ on $\mathcal{P}_+$ we define by
$$
\delta_x \bowtie \delta_y = \frac{1}{2} \delta_{x+y} + \frac{1}{2} \delta_{|x-y|}, \quad x,y \geqslant 0.
$$
This distribution can be considered as the limit of $\delta_x \otimes_{\omega_s}  \delta_y$ for $s \searrow -\frac{1}{2}$.
\end{Example}
\vspace{-5mm}
\begin{Example}
$\alpha$-stable convolution $\ast_{\alpha}$ for $\alpha >0$ is given for $x,y \geqslant 0$ by
$$ 
\delta_x \ast_{\alpha} \delta_y = \delta_{g_{\alpha}(x,y)}, \quad \hbox{ where } \quad g_{\alpha}(x,y) = (x^{\alpha} +y^{\alpha} )^{1/{\alpha}}.
$$
\end{Example}
\vspace{-8mm}
\begin{Example}
The Kendall generalized convolution  $\vartriangle_{\alpha}$ on $\mathcal{P}_+$, $\alpha>0$, is defined by
$$
\delta_x \vartriangle_{\alpha} \delta_1 := \bigl(1- x^{\alpha} \bigr) \delta_1 + x^{\alpha} \pi_{2\alpha},\quad x\in[0,1],
$$
where $\pi_{\beta}$ is the Pareto distribution with the density function $f_{\beta}(t) = \beta t^{-\beta-1}$ on the set $[1,\infty)$ for some $\beta>0$.
\end{Example}
\vspace{-5mm}
\begin{Example}
The $\max$-convolution is simply defined by 
$$
\delta_x \triangledown \delta_y = \delta_{x \vee y}.
$$
This distribution can be considered as the limit of 
$\delta_x \vartriangle_{\alpha} \delta_y$ for $\alpha \rightarrow \infty$. This is the reason why we call the Kendall convolution by the extremal Kendall convolution.
\end{Example}
\vspace{-5mm}
\begin{Example}
The Kucharczak convolution $\delta_x \circ_1 \delta_y$ for $x,y  \geqslant 0$ defined in \cite{Urbanik93}, Example 2.4, is a measure absolutely continuous with respect to the Lebesgue measure and for $a\in (0,1]$, $r >0$, given by
$$
\delta_x \circ_1 \delta_y (dt) = \frac{r x^a y^a}{\Gamma(a) \Gamma(1-a)} \, \frac{t^{r-ar-1} (2t^r -x^r - y^r) \mathbf{1}_{[g_{ar}(x,y), \infty)}(t)}{(t^r -x^r -y^r)^a (t^r -x^r) (t^r - y^r)} \, dt.
$$
\end{Example}
\vspace{-5mm}
\begin{Example}\label{KUr}
The Kucharczak-Urbanik  generalized convolution  defined in \cite{KucharczakUrbanik, Urbanik93} for $\alpha >0$ and $n \in \mathbb{N}$ is uniquely determined by
\begin{equation}\label{KUdef}
\delta_x \vartriangle_{\alpha,n} \delta_1 (ds)  
:= (1 - x^{\alpha})^{n} \delta_1(ds) + \sum_{k=1}^n {{n}\choose {k}} x^{\alpha k} (1 - x^{\alpha})^{n-k} \mu_{k,n} (ds)
\end{equation}
for $x\in[0,1]$, where for $k=1,\ldots,n$  the probability measures $\mu_{k,n}$, are defined by their density functions $f_{k,n}$:
\begin{equation}\label{efkaen}
f_{k,n}(s) = \alpha k {{n+k}\choose n} s^{-\alpha(n+1) - 1} \left( 1 - s^{-\alpha} \right)^{k-1},\quad s>1.
\end{equation}
\end{Example}
\vspace{-5mm}
\begin{Example}
A family of non-regular generalized convolutions $\diamondsuit_{p,\alpha}$, $p\in[0,1]$, $\alpha >0$, was introduced by K. Urbanik in \cite{Urbanik93} initially for $\alpha =1$. This family interpolates between two boundary cases: the max-convolution for $p=0$ and the Kendall convolution for $p=1$.
The $\diamondsuit_{p,\alpha}$-convolution $\delta_x \diamondsuit_{p,\alpha} \delta_1$, $x \in [0,1]$, is defined  for $p \neq \frac{1}{2}$ by 
$$
\delta_x \diamondsuit_{p,\alpha} \delta_1 (ds)= (1-px^{\alpha})\, \delta_1(ds) +px^{\alpha} \,\tfrac{\alpha}{2p-1} \tfrac{2p - s^{q}}{s^{2\alpha+1}} \mathbf{1}_{[1,\infty)}(s) ds,\quad x\in[0,1],
$$
where $q = \frac{\alpha(1-2p)}{(1-p)}$. By continuity, for $p\to {1/2}$ we have
$$
\delta_x \diamondsuit_{1/2,\alpha} \delta_1 (ds)= \bigl(1-\tfrac{x^{\alpha}}{2} \bigr)\, \delta_1(ds) + \tfrac{x^{\alpha}}{2}\, \tfrac{\alpha( 1 + 2 \ln s)}{s^{2\alpha+1}}\mathbf{1}_{[1,\infty)}(s) ds.
$$
\end{Example}
\begin{Example}
In  \cite{Mej18l-c} one can find the  description of the regular generalized convolutions called the Kendall-type convolutions. Their probability kernels  are the following:
$$
\varphi_{c,\alpha,p}(t) = \left( 1 - (1+c) t^{\alpha} + ct^{\alpha p} \right) \mathbf{1}_{[0,1]}(t),
$$
where $p \geqslant 2$, $\alpha >0$ and one of the following conditions holds\\
{\bf 1)} $c= (p-1)^{-1}$, \\
{\bf 2)} $c = (p^2 -1)^{-1}$, \\
{\bf 3)} $c = \frac{1}{2} (2-p) (p-1)^{-1}$, \\
{\bf 4)} $c = \frac{1}{2} (p-1)^{-1}$, \\
{\bf 5)} $c \in \bigl( (p^2 -1)^{-1}, \frac{1}{2} (p-1)^{-1} \bigr)$ and none of the previous cases holds.

For other parameters $p, c, \alpha$ none of the functions $\varphi_{c,\alpha,p}$ can be a probability kernel of a regular generalized convolution. Such convolutions are given by
$$
\delta_x \vartriangle_{c,\alpha,p} \delta_1 = \varphi_{c,\alpha,p}(x) \, \delta_1 + x^{\alpha p} \, \lambda_1 + (c+1)(x^{\alpha} -x^{\alpha p})\, \lambda_2,
$$
for properly chosen probability measures $\lambda_1, \lambda_2$ supported in $[1,\infty)$. For details, in particular for the explicit densities and cumulative distribution functions of the measures $\lambda_1, \lambda_2$,  see \cite{Mej18l-c}.
\end{Example}
\section{The Kendall convolution by the corresponding Williamson transform}

Let $m_0$ denote the sum of $\delta_0$ and the Lebesgue measure on $[0,\infty)$. By Theorem 4.1 and Corollary 4.4 in \cite{Urbanik86} we know that the generalized convolution can be defined uniquely by its generalized characteristic function treated as an integral transform. Such approach is described by the next definition. Let us remind here that the $L_1(m_0)$-topology of $L_{\infty}(m_0)$-space means that $f_n \rightarrow f$ for functions $f_n, f \in L_{\infty}(m_0)$ if $\int f_n(x) g(x) m_0(dx) \rightarrow \int f(x) g(x) m_0(dx)$ for all $g \in L_1(m_0)$, or, equivalently if 
$$
\int f_n(yx) g(x) m_0(dx) \rightarrow \int f(yx) g(x) m_0(dx)
$$ 
for all $g \in L_1(m_0)$ and all $y\in [0,\infty)$.

\begin{defn}
We say that the Borel function $\varphi \colon [0,\infty) \rightarrow \mathbb{R}$, $|\varphi(t) | \leqslant \varphi(0) = 1$, defines a $\varphi$-generalized convolution on $\mathcal{P}_+$ if \\
{\bf (i)}  the integral transform
$$
\widehat{\lambda} (t) := \int_0^{\infty} \varphi(tx) \lambda(dx), \quad \lambda \in \mathcal{P}_+,
$$
separates points in $\mathcal{P}_+$, i.e. $\widehat{\lambda}= \widehat{\mu}$ implies that $\lambda = \mu$, \\
{\bf (ii)} the weak convergence $\lambda_n \rightarrow \lambda$ is equivalent to the convergence $\widehat{\lambda_n} \rightarrow \widehat{\lambda}$ in the $L_1(m_0)$-topology of $L_{\infty}(m_0)$, \\
{\bf (iii)} for every $x,y \geqslant 0$ there exists a measure $\mu \in \mathcal{P}_+$, such that the following equality, called the product formula for the function $\varphi$, holds
\begin{equation} \label{productformula}
\forall \, x,y \geqslant 0 \, \exists \, \mu \in \mathcal{P}_+ \qquad \varphi(xt) \, \varphi(yt) = \int_0^{\infty}\! \varphi(st)\, \mu(ds).
\end{equation}
The corresponding $\varphi$-generalized convolution for such function $\varphi$ and measure $\mu$ described in {\rm (iii)} is defined by 
$$
\forall \, x,y \in \mathbb{R}_+ \quad  \delta_x \diamondsuit_{\varphi} \delta_y :=  \mu. 
$$
\end{defn} 

\noindent
\begin{rem}
It is easy to notice that the operation $\diamondsuit_{\varphi}$ defined for the point mass measures by the Definition 4 satisfies all  conditions of Definition 1 thus it is a generalized convolution in the Urbanik sense.  Continuity of the function $\varphi$ is equivalent with the regularity of convolution $\diamondsuit_{\varphi}$.
\end{rem}

In all the following examples, except  the Example 4.10, we see that the known generalized convolution $\diamond$  is $\varphi$-generalized convolution for the function $\varphi$ being the probability kernel for $\diamond$. In Example 4.10 we describe Whittaker $W_{\alpha,\nu}$-generalized convolution, based on a little changed product formula. This convolution does not satisfy the condition ($iii$) of the Urbanik definition of generalized convolution, however  methods used in studying the properties of one convolutions can be helpful in studying the properties of others. \\

{\bf Example 4.0.} The characteristic function of the variable $\theta_s$, $s > - \frac{1}{2}$ is given by the following formula (for the proof see e.g. \cite{King})
$$
\Phi_{s} (t) := \int_{-1}^1 e^{itx} f_s(x) dx = \Gamma(s+1) \Bigl( \frac{t}{2} \Bigr)^s J_s(t),
$$
where $J_s$ is the Bessel function of the first kind with the index $s$ and
$$
J_s(t) = \sum_{m=0}^{\infty} \frac{(-1)^m }{m! \Gamma(m+1 +s)} \ \Bigl( \frac{t}{2} \Bigr)^{2m -s}.
$$
We recognize here $\varphi = \Phi_s$, $\widehat{\lambda}(t) = \int_0^{\infty} \Phi_s(tx) \lambda(dx)$. The definition of the Kingman-Bessel convolution $\otimes_{\omega_s}$ follows now from the Gegenbauer's Formula (see e.g. \cite{tich}, Chapter 8.19), which is the product formula (\ref{productformula}) for this case: 
$$
\Phi_s (xt) \, \Phi_s (yt) = \int_0^{\infty} \Phi_s (rt) \, r_{s}(x,y,r) \, dr.
$$
Here for $x,y>0$ the function $r_s(x,y,r)$ as a function on $r$, is the density of the random variable $\sqrt{x^2 + y^2 + 2xy \theta_s}$ and it is equal to
$$
r_s(x,y,r) = \frac{\Gamma(s+1)}{\sqrt{\pi} \Gamma( s + {1/2}) } \, \frac{2^{1-2s} (xy)^{-2s} \, r \, \mathbf{1}_{(|x-y|, x+y)}(r) }{\bigl[ (r^2 - (x-y)^2) ((x+y)^2 - r^2 ) \bigr]^{-s+ \frac{1}{2}} }.
$$
\begin{Example}
For the classical convolution we have $\varphi(t) = e^{-t}$ and the integral transform $\lambda \rightarrow \widehat{\lambda}$ is the classical Laplace transform. The product formula (\ref{productformula}) follows from the fact that the Laplace transform of the convolution of two measures is equal to the product of their Laplace transforms.
\end{Example}
\begin{Example}
For the symmetric convolution $\bowtie$ we have $\varphi(t) = \cos(xt)$ and equation (6) follows from the elementary formula:
$$
\cos(xt) \cos(yt) = \frac{1}{2} \cos((x+y)t) + \frac{1}{2} \cos((x-y)t).
$$
\end{Example}
\begin{Example}
For the $\alpha$-stable convolution $\ast_{\alpha}$ we have $\varphi(t) = e^{-t^{\alpha}}$. This means that $\widehat{\lambda}$ is simply a modified Laplace transform.
\end{Example}
\begin{Example}
Recall that for $\alpha >0$ and  a non-negative, $\sigma$-finite on $[0,\infty)$ (finite on compact sets)  measure $\lambda$ on $\mathbb{R}_+$ the Williamson transform $\mathcal{W}_{\alpha}\lambda$ is  defined by 
$$
\mathcal{W}_{\alpha} \lambda(t):= \int_0^{\infty} \bigl( 1 - t^{\alpha} x^{\alpha} \bigr)_+ \lambda(dx),
$$
where $a_+=\max\{a,0\}$. The product formula (\ref{productformula}) for the Williamson transform is the following: 
$$
(1-t^{\alpha} x^{\alpha})_+ (1-t^{\alpha} y^{\alpha})_+ = \int_0^{\infty}\!\! (1-t^{\alpha} s^{\alpha} )_{+} (\delta_x \vartriangle_{\alpha} \delta_y) (ds), \qquad x,y \geqslant 0.
$$
This formula was introduced for studying of the Kendall $\vartriangle_{\alpha}$ convolution
in \cite{Urbanik84}, thus $\vartriangle_{\alpha}$-generalized characteristic function is given by:
\begin{equation}\label{wil}
\Phi^{\vartriangle_{\alpha}}_{\lambda}(t) := \mathcal{W}_{\alpha} \lambda(t)= \int_0^{\infty} \bigl( 1 - t^{\alpha} s^{\alpha} \bigr)_+ \lambda(ds).
\end{equation}

The Williamson integral transform for $\alpha =1$ was introduced when studying $n$-times monotonic functions, i.e. functions $f$ on $[0,\infty)$ such that $(-1)^{\ell} f^{(\ell)}(r)$ is non-negative, non-decreasing and convex for $\ell = 0,1,...,n-1$. R.E. Williamson showed (see \cite{Williamson} Th. 1 and 2) that $f$ is $n$-times monotonic function on $(0,\infty)$ iff $f(t) = \int_0^{\infty} (1-tx)_+^{n-1} \gamma(dx)$,
for some non-negative, $\sigma$-finite measure $\gamma$ on $[0,\infty)$.

Actually, the original Williamson transform and its modifications
$\gamma \longrightarrow \int_0^{\infty} \bigl( 1- t^{\alpha} x^{\alpha} \bigr)_+^{d-1}\, \gamma(dx)$,
for some $\alpha, d >0$, are applied in many different areas of mathematics including actuarial science (see e.g. \cite{Castaner, Lefevre}) and dependence modeling by copulas (\cite{GenNesRiv, Larsson,McNeilNesholva1, McNeilNesholva2}). 

Note that it is easy to retrieve the measure knowing its Williamson transform. This makes the proof of the fact that the Williamson transform uniquely determines the measure much simpler than that for the Fourier or Laplace transforms.  To see this we integrate by parts the right hand side of \eqref{wil} and we obtain
$$
\Phi_{\lambda}^{\vartriangle_{\alpha}} (t) = \alpha t^{\alpha} \int_0^{1/t} x^{\alpha -1} F(x)\, dx,
$$
where $F$ is the cumulative distribution function for $\lambda$. Now,  with the notation $G(t) = \Phi_{\lambda}^{\vartriangle_{\alpha}}(1/t)$, we obtain
\begin{equation}\label{F_lambda}
t^{\alpha} G(t) = \alpha \int_0^t x^{\alpha -1} F(x) dx, \quad \hbox{ thus } \quad F(t) = G(t) + \alpha^{-1} t^{-1} G'(t),
\end{equation}
at each continuity point of the function $F$. Consequently, $\Phi_{\lambda_1}^{\vartriangle_{\alpha}} (t) = \Phi_{\lambda_2}^{\vartriangle_{\alpha}} (t)$ implies that $\lambda_1 = \lambda_2$.
Since  $\Phi_{\lambda}^{\vartriangle_{\alpha}} (t)$ is the generalized characteristic function for the Kendall convolution  we know that for $\lambda_1, \lambda_2 \in \mathcal{P}_+$  \begin{equation}\label{will}
\Phi_{\lambda_1 \vartriangle_{\alpha} \lambda_2}^{\vartriangle_{\alpha}}(t) = \Phi_{\lambda_1}^{\vartriangle_{\alpha}} (t) \, \Phi_{\lambda_2}^{\vartriangle_{\alpha}} (t), \qquad \qquad t \geqslant 0.
\end{equation} 
The cumulative distribution function of the Kendall convolution of two measures can also be be easily expressed:

\begin{thm}
For every $\lambda_1, \lambda_2 \in \mathcal{P}_+$. The measure $\lambda= \lambda_1\vartriangle_{\alpha} \lambda_2$  if and only if $F_{\lambda}$, the cumulative  distribution function of $\lambda$, is given by
$$
F_{\lambda}(t) = G_1 (t) F_2 (t) + G_2 (t) F_1(t) - G_1(t) G_2(t),
$$
where $F_i$ is the cumulative distribution function of $\lambda_i$, $G_i(t) = \Phi_{\lambda_i}^{\vartriangle_{\alpha}}(1/t)$, $i=1,2$, $F_{\lambda}$ is the cumulative distribution function of $\lambda$ and $G_{\lambda} = \Phi_{\lambda}^{\vartriangle_{\alpha}}(1/t)$.
\end{thm}

\noindent
{\bf Proof.}
Assume that $\lambda = \lambda_1 \vartriangle_{\alpha} \lambda_2$.
By the formula expressing the cumulative distribution function by the Williamson transform and the equality 
$\widehat{\lambda_1 \vartriangle_{\alpha} \lambda_2} = \widehat{\lambda_1} \widehat{\lambda_2}$ we have that $G_{\lambda}(t) = G_1(t) G_2(t)$, $t\geqslant 0$, and then, by \eqref{F_lambda}
\begin{eqnarray*}
F_{\lambda}(t) & = & G_{\lambda}(t)) + \alpha^{-1} t^{-1} G_{\lambda}'(t) \\
& = & G_{1}(t) G_{2}(t) + \alpha^{-1} t^{-1} G_{1}'(t) G_{2}(t) + \alpha^{-1} t^{-1} G_{1}(t) G_{2}'(t) \\
& = & G_1 (t) F_2 (t) + G_2 (t) F_1(t) - G_1 (t) G_2(t).
\end{eqnarray*}
Assume now that the cumulative distribution function $F_{\lambda}$ can be written by the desired formula. Since $F_i(t) = G_i(t) + \alpha^{-1} t^{-1} G_i'(t)$, $i=1,2$, then
\begin{eqnarray*}
F_{\lambda}(t) & = & G_1 (t) F_2 (t) + G_2 (t) F_1(t) - G_1(t) G_2(t) \\
& = & G_1(t) G_2(t) + \alpha^{-1} t^{-1} \bigl(G_{1} (t) G_2 (t) \bigr)'.
\end{eqnarray*}
By the uniqueness of the Williamson transform we see that the generalized characteristic function of $\lambda$ is equal to $G_1(t^{-1}) G_2(t^{-1})$, $t \geqslant 0$ which is the generalized characteristic function of $\lambda_1 \vartriangle_{\alpha} \lambda_2$. \qed 
\end{Example}

\begin{Example}
For the $\max$-convolution we have $\varphi(t) = \mathbf{1}_{[0,1]}(t)$. This function is not continuous, thus the corresponding convolution $\triangledown$ is not regular, but the inversion formula is also equally easy to obtain:
$$
\widehat{\lambda}(t) = \int_0^{\infty} \mathbf{1}_{[0,1]}(tx) \lambda(dx) = \int_0^{t^{-1}}  \lambda(dx)=  F_{\lambda}(t^{-1}),
$$
thus $F_{\lambda}(t) = \widehat{\lambda}(t^{-1})$ for all continuity points of the cumulative distribution function $F_{\lambda}$.
\end{Example}

\begin{Example}
For $a\in (0,1]$, $r>0$, the Kucharczak generalized convolution $\circ_1$ can be defined by the product formula (\ref{productformula})  applied to its probability kernel:
$$
\Omega(t) = \frac{\Gamma(a,t^r)}{\Gamma(a)} = \frac{1}{\Gamma(a)} \int_{t^r}^{\infty} x^{a-1} e^{-x} dx, \quad \quad t>0.
$$
This means that the measure $\mu = \delta_x \circ_1 \delta_y$ is defined as a solution of the following integral equation:
$$
\frac{1}{\Gamma(a)^2} \int_{t^r x^r}^{\infty} \! s^{a-1} e^{-s} ds \, \int_{t^r y^r}^{\infty} \! u^{a-1} e^{-u} du = \int_0^{\infty} \!\!\! \frac{1}{\Gamma(a)} \int_{t^r s^r}^{\infty} \! u^{a-1} e^{-u} du \,\mu(ds).
$$
\end{Example}
\begin{Example}
The Kucharczak-Urbanik convolution $\vartriangle_{\alpha,n}$ can be defined by equation (\ref{productformula}) for $\varphi(t) := (1-t^{\alpha})_+^n$. To see this note that for any $x\in [0,1]$ and $t\geqslant 0$ we have
$$
(1-t^{\alpha} x^{\alpha})_+^n (1-t^{\alpha})_+^n  = \sum_{k=0}^{n} {{n}\choose{k}} x^{\alpha k}  (1 - x^{\alpha})^{n-k} (1-t^{\alpha})_+^{n+k}.
$$
It remains to show that for any integer $k\ge 1$
$$
(1-t^{\alpha})_+^{n+k} = \int_0^{\infty}\!\!  (1-t^{\alpha} s^{\alpha})_+\, f_{k,n}(s)\, ds, 
$$
where the density functions $f_{k,n}$, $k=1,\ldots,n$, $n \in \mathbb{N}$, are described in Example 3.7. This equality we can obtain by a simple induction argument (with respect to $k$), where the first step  of induction is based on the following property of the Pareto distribution:
$$
\int_0^{\infty}\!(1- s^{\alpha} t^{\alpha})_+^n \, \pi_{\alpha (n+1)}(ds) = (1-  t^{\alpha})_+^{n+1}.
$$
The final conclusion is a consequence of the uniqueness of probability kernel (up to a scale coefficient) of every generalized convolution (for the proofs see \cite{MisiewiczVolodia1,Urbanik84}). 

The inversion formula for the integral transform $\lambda \rightarrow \widehat{\lambda}$ with the kernel $\Omega_{\alpha, n}$ can be obtained using the same methods as inverting the Williamson transform, but the level of difficulty increases with the increase of $n$ - for the detailed proof see \cite{McNeilNesholva1, McNeilNesholva2} .
\end{Example}
\begin{Example}
The $\diamondsuit_{p,\alpha}$ generalized convolution, $\alpha >0$, $p\in [0,1]$,  can be defined by equation (\ref{productformula}) for the  probability kernel
$$
\Omega_{\diamondsuit_{\alpha, p}}(t) = (1-pt^{\alpha}){\bf 1}_{[0,1]}(t),\quad t\ge 0
$$
This function, except for the Kendall case $p=1$, is not continuous thus the generalized convolution $\diamondsuit_{p,\alpha}$ is not regular. 
\end{Example}
\begin{Example}
The Kendall-type generalized convolutions $\vartriangle_{c,\alpha,p}$ were found by considering such parameters $c,\alpha, p$ for which the function $\varphi_{c,\alpha,p}(t) = \left( 1 - (1+c) t^{\alpha} + ct^{\alpha p} \right) \mathbf{1}_{[0,1]}(t)$ can play the role of probability kernel of some generalized convolution. In particular we choose such $c,\alpha, p$ that for all $x,y>0$ the measure $\mu$ (depending on $x$ and $y$) which satisfies the equality 
$$
\varphi_{c,\alpha,p}(xt) \, \varphi_{c,\alpha,p}(yt) = \int_0^{\infty} \! \varphi_{c,\alpha,p}(ts) \mu(ds)
$$
is a probability measure. 
\end{Example}

\subsection{Generalized convolutions in Harmonic Analysis}.

\vspace{3mm}

The version of equation (\ref{productformula}) appearing in the theory of special functions and harmonic analysis is called a product formula or a multiplication formula for the family $\{ \chi_{_{\lambda}}\}_{\lambda \in \Lambda}$  of continuous functions on $I \subset \mathbb{R}$:
$$
\chi_{_{\lambda}}(x) \, \chi_{_{\lambda}}(y) = \int_{I} \chi_{_{\lambda}}(s) K(x,y,s) \, ds, \quad \lambda \in \Lambda, \leqno{(6')}
$$
where the kernel $K(x,y,s)$ does not depend on $\lambda$ and $\Lambda$ is some indexing  set. Such product formulas are the key ingredient for definitions of {\it generalized translation} and {\it generalized convolution operators} which have been  introduced by J. Delsarte \cite{Delsarte} and B. Levitan \cite{Levitan} in the theory of special functions and harmonic analysis. For details and examples see \cite{Berezansky, Connett}.

For the generalized convolutions on $\mathcal{P}_+$ introduced by K. Urbanik in the probability theory we have
$$
\bigl\{ \chi_{_{\lambda}}(\cdot) \colon \lambda \in I \bigr\} = \bigl\{ \Omega(t\, \cdot) \colon t \geqslant 0 \bigr\},
$$
where $\Omega \colon \mathbb{R}_+ \to \mathbb{R}_+$ is the probability kernel for the considered generalized convolution. In the definition of  J. Delsarte \cite{Delsarte} and B. Levitan \cite{Levitan} the set $I$ in the family $\{ \chi_{_{\lambda}}\}_{\lambda \in \Lambda}$ is some indexing set and the equality $\chi_{_{\lambda}}(x) = \chi_{_{1}}(\lambda x)$ does not have to hold, but the family
$$
\left\{ \int_{I} \chi_{_{\lambda}}(s) K(x,y,s) \, ds  \colon \lambda \in \Lambda \right\}
$$
will identify the kernel $K(x,y,s)$ uniquely up to a set of Lebesgue measure zero for each choice of $x,y \in I$. \\

{\bf Example 4.10.} In \cite{Sousa1(19)} the authors proved the product formula for the index Whittaker transform and defined the corresponding generalized convolution operator. By the index Whittaker transform we understand here the integral transform $\mathcal{P}_+ \ni \mu \rightarrow (W_{\alpha,\nu} \mu)(t)$ given by
$$
\widehat{\mu}(\lambda) := (W_{\alpha} \mu) (\lambda) := \int_0^{\infty} \! W_{\alpha, \Delta_{\lambda}}(x)  \mu(dx), \qquad \lambda \geqslant 0,
$$
where $\alpha < \frac{1}{2}$ is a parameter, $\Delta_{\lambda} = \sqrt{(\frac{1}{2} - \alpha)^2 - \lambda}$ and $W_{\alpha, \nu}$ is the Whittaker function
$$
W_{\alpha, \nu}(x) = \frac{e^{-\frac{x}{2}} x^{\alpha}}{\Gamma( \frac{1}{2} - \alpha + \nu)} \int_0^{\infty}\! e^{-s} s^{- \frac{1}{2} -\alpha + \nu} \Bigl( 1 + \frac{s}{x}\Bigr)^{-\frac{1}{2} + \alpha + \nu} ds,
$$
for $\mathfrak{R}e \, x >0, \mathfrak{R}e \, \alpha < \frac{1}{2} + \mathfrak{R}e \,\nu$. Equivalently the Whittaker function is defined as the solution of Whittaker's differential equation:
$$
\frac{d^2 u}{d x^2} + \Bigl( - \frac{1}{4} + \frac{\alpha}{x} + \frac{{1/4} - \nu^2}{x^2} \Bigr) u = 0
$$
uniquely determined by the property $W_{\alpha, \nu}(x) \sim x^{\alpha} e^{-{x/2}}$ for $x \rightarrow 0$. 

The index Whittaker transform $\mu \rightarrow \widehat{\mu}$ has the following properties of the generalized characteristic function (see Prop. 4.4 in \cite{Sousa2(19)}): \\

($i$) $\widehat{\mu}$ is uniformly continuous on $[0,\infty)$. Moreover, for any indexing set $J$ if the family $\{ \mu_j|_{(0,\infty)} \colon j \in J \}$ is tight, then $\{ \widehat{\mu_j} \colon j \in J\}$ is uniformly equicontinuous; \\

($ii$) $\widehat{\mu}$ uniquely determines $\mu \in \mathcal{P}_+$; \\

($iii$) if $\mu_n, \mu \in \mathcal{P}_+$, $n \in \mathbb{N}$, and $\mu_n \Rightarrow \mu$ then $\widehat{\mu_n} \rightarrow \widehat{\mu}$ uniformly on compact sets; \\

($iv$) if $\mu_n \in \mathcal{P}_+$, $n \in \mathbb{N}$ and $\widehat{\mu_n}(\lambda) \rightarrow f(\lambda)$ pointwise in $\lambda \geqslant 0$ for some real function $f$, continuous in a neighbourhood of zero then there exists $\mu \in \mathcal{P}_+$ such that $f = \widehat{\mu}$. \\

The product formula for the Whittaker function of the second kind is the following (see Th. 3.1  in \cite{Sousa1(19)}):
\begin{equation} \label{Whittakereq}
W_{\alpha,\nu}(x) \, W_{\alpha,\nu}(y) = \int_0^{\infty}\! W_{\alpha,\nu}(s) K_{\alpha} (x,y,s) ds,
\end{equation}
where
\begin{eqnarray*}
\lefteqn{K_{\alpha} (x,y,s) := \frac{(xy)^{2\alpha-1}}{\sqrt{2\pi} s^{2\alpha}} \, \times} \\
&& \hspace{-5mm} \exp\left\{ \frac{1}{2x^2}+ \frac{1}{2y^2} - \frac{1}{2s^2} - \Bigl(\frac{x^2 + y^2 + s^2}{4xys}\Bigr)^2 \right\} D_{2\alpha} \left[\frac{x^2 + y^2 + s^2}{{2xys}} \right]
\end{eqnarray*}
and $D_{\mu}(s)$ is the parabolic cylinder function for $s>0$, $\mathfrak{R}e \, \mu <1$:
$$
D_{\mu}(s) = \frac{s^{\mu} e^{-{s^2}/4}}{\Gamma( \frac{1}{2} (1-\mu))} \int_0^{\infty} t^{\frac{1}{2}(1+\mu)} \Bigl( 1 + \frac{2t}{s^2}\Bigr)^{{\mu}/2}\, e^{-t}  ds.
$$
The equation (\ref{Whittakereq}) holds for all $\nu$ for which the function $W_{\alpha, \nu}$ can be defined, but considering generalized characteristic function in the sense of Delsarte and Levitan we will assume that $\nu = \Delta_{\lambda}$.  
By Theorem 4.6 in \cite{Sousa2(19)} we have $\int_0^{\infty} K_{\alpha} (x,y,s) ds = 1$ for all $x,y >0$. Consequently we have that the product formula (\ref{Whittakereq}) for the Whittaker function defines a generalized convolution $\maltese$ in the sense of Delsarte and Levitan:
$$
\delta_x \,\maltese \, \delta_y (ds) = K_{\alpha}(x,y,s) \mathbf{1}_{(0,\infty)} (s) \, ds.
$$
This proposal does not guarantee that $\maltese$ is a generalized convolution in the Urbanik's sense. In particular, we do not know if  conditions ($iii$) or ($v$) of Definition 1 hold.

\section{The Kendall convolution as a weak generalized convolution}
Let us remind that the measure $\nu \in \mathcal{P}(\mathbb{E})$ is stable if for all $a,b \geqslant 0$ there exists non-random  $\mathbf{d}(a,b) \in \mathbb{E}$ such that
$$
T_a \nu \ast T_b \nu = T_{c(a,b)} \nu \ast \delta_{\mathbf{d}(a,b)},
$$
where $c(a,b)^{\alpha} = a^{\alpha} + b^{\alpha}$ for some $\alpha \in (0,2]$. If $\mathbf{d}(a,b) \equiv 0$ then the measure $\nu$ is called strictly stable.
The complete characterization of both stable and strictly stable distributions is known and given e.g. in \cite{ST}.

Similarly we define  weakly stable distributions, which are measures on an arbitrary separable Banach space $\mathbb E$ (with the Borel $\sigma$-algebra):
\begin{defn}\label{ws} 
We say that a measure $\mu \in \mathcal{P}(\mathbb{E})$  is weakly stable if
$$
\forall\, a,b \in \mathbb{R} \,\, \exists \, \lambda = \lambda_{a,b} \in \mathcal{P}: \qquad  T_a \mu \ast T_b \mu = \lambda \circ \mu,
$$
where $\ast$ denotes the classical convolution and $\mathcal{L}(X) \circ \mathcal{L}(\theta) = \mathcal{L}(X\theta)$ if the random elements $X$ and $\theta$ are independent. 
\end{defn}

It is known (see \cite{MOU}) that $\mu$  is weakly stable if and only if  
$$
\forall\, \lambda_1, \lambda_2 \in \mathcal{P} \,\, \exists \, \lambda = \lambda_{a,b} \in \mathcal{P} \quad  (\lambda_1 \circ \mu) \ast (\lambda_2 \circ \mu) = \lambda \circ \mu.\eqno{(\ast)}
$$
This property is the base for defining weak generalized convolution:
\begin{defn}
Let $\mu$ be a weakly stable measure on a separable Banach space $\mathbb{E}$. The binary operation $\otimes_{\mu} \colon \mathcal{P}_+^2 \rightarrow \mathcal{P}_+$, called a $\mu$-weak generalized convolution, is defined as follows: for any $\lambda_1,\lambda_2\in\mathcal{P}_+$
$$
\lambda_1 \otimes_{\mu} \lambda_2 = \lambda\in\mathcal{P}_+ \quad\Longleftrightarrow\quad \bigl(\lambda_1 \circ \mu \bigr) \ast \bigl(\lambda_2 \circ \mu \bigr) = \lambda \circ \mu.
$$
The generalized convolution $\diamond$ is called a weak generalized convolution if there exists a weakly stable measure $\mu$ such that $\diamond = \otimes_{\mu}$.
\end{defn}
All known weakly stable measures are symmetric, i.e. satisfying the property $\mu(A) = \mu(-A)$ for every Borel set $A \in \mathcal{B}(\mathbb{E})$. Moreover if $\mu$ on $\mathbb{E}$ is weakly stable, then for every subspace $\mathbb{E}_1 \subset \mathbb{E}$ and every linear operator $Q \colon \mathbb{E} \rightarrow \mathbb{E}_1$ the measure $\mu_Q$ on $\mathbb{E}_1$ defined by
$$
\forall \, A \in \mathcal{B}(\mathbb{E}) \quad \mu_Q(A) := \mu (Q^{-1}(A))
$$
is also weakly stable and  both $\mu$ and $\mu_Q$  define the same weak generalized convolution on $\mathcal{P}_+$. For these reasons in defining weak generalized convolutions we will restrict our attention to weakly stable measures $\mu \in \mathcal{P}_s$ (symmetric measures on $\mathbb{R}$). \\
\begin{rem} {\rm
Let $\widehat{\mu}$ be the characteristic function of the weakly stable measure $\mu \in\mathcal{P}_s$. By the weak stability condition ($\ast$) and Definition 6  written in the language of characteristic functions we have that there exists a measure $\lambda= \lambda_1 \otimes_{\mu} \lambda_2$ such that
$$
\int_0^{\infty} \widehat{\mu }(st) \lambda_1(ds) \int_0^{\infty} \widehat{\mu}(st ) \lambda_2(ds) = \int_0^{\infty} \widehat{\mu}(st) \lambda(ds).
$$
In this case the probability kernel of the generalized convolution $\otimes_{\mu}$ is $\varphi(t) = \widehat{\mu}(t) =  \int_{\mathbb{R}}\, \cos(tx) \mu(dx)$ considered as a function on $\mathbb{R}_+$. Finally we see that the generalized convolution $\diamond$ with the probability kernel $\varphi$ is weak iff the function $\varphi(|t|)$, $t \in \mathbb{R}$, is a characteristic function of some measure $\mu \in \mathcal{P}_s$, $\diamond = \otimes_{\mu}$ and
\begin{equation}\label{weakdef}
\forall \, a,b,t \in \mathbb{R}_+ \quad \varphi(at) \varphi(bt) = \int_0^{\infty} \varphi(st) \, \delta_a \otimes_{\mu} \delta_b\,  (ds).
\end{equation} }
\end{rem}

\begin{thm}
The Kendall convolution $\vartriangle_{\alpha}$ is a  weak generalized convolution if  $\alpha \in(0,1]$. The corresponding weakly stable measure $\mu:=\mu_{\alpha}\in\mathcal{P}_s$ is defined by the density function
$$
g_{\alpha}(t) = \frac{2\alpha}{\pi}\, |t|^{-\alpha -1} \int_0^{|t|} x^{\alpha -1} \sin{x} \, dx,\quad t\in\mathbb{R}\setminus\{0\}.
$$
\end{thm}
\noindent{\bf Proof.}
Since we already know that the probability kernel for the Kendall convolution is  $\Omega_{\vartriangle_{\alpha}}(t) = (1- t^{\alpha}) \mathbf{1}_{[0,1]}(t)$  we only need to:  
\begin{namelist}{lll}
    \item{a)} show that the function  $g(t) := \Omega_{\vartriangle_{\alpha}}(|t|)$ is a characteristic function of some probability $\mu$ if  $\alpha \in (0,1]$; 
    \item{b)} identify $\widehat{g}_{\alpha}$ as the density of $\mu$. 
\end{namelist}
Indeed, if a) and b) hold then $\widehat{\mu} (t) = \Omega_{\vartriangle_{\alpha}}(|t|)$ thus, by equality (\ref{weakdef}) we see that $\mu$ is weakly stable and defines the convolution $\vartriangle_{\alpha}$.  \\
To see that a) holds true, note that for $t>0$
$$
\widehat{g}'(t) = - \alpha t^{\alpha-1} < 0, \quad \hbox{ and } \quad \widehat{g}''(t) = \alpha (1-\alpha) t^{\alpha -2} > 0.
$$
Consequently, by the Polya Theorem, it follows that $\widehat{g}$ is indeed the characteristic function of a symmetric probability measure $\mu$. 

To see that b) holds true we use  the inverse Fourier transform for integrable characteristic function to obtain the density function  of $\mu$:
\begin{align*}
\frac{1}{2\pi}\int_{\mathbb{R}} g(x) e^{-itx} \, dx & = \!\int_{\mathbb{R}} \bigl( 1 - |x|^{\alpha} \bigr)_+ e^{-itx} \, dx  \\
&= \frac{\alpha}{\pi} \, |t|^{-\alpha -1} \! \int_0^{|t|} x^{\alpha -1} \sin{x} \, dx=g_{\alpha}(t).  \qquad \qquad{\Box}
\end{align*} 

\begin{thm}
Assume that the Kendall convolution $\vartriangle_{\alpha}$ is weakly stable. Then $\alpha \in (0,2]$.
\end{thm}

\noindent{\bf Proof.}
The Kendall convolution $\vartriangle_{\alpha}$ is weakly stable iff the function $(1 - |t|^{\alpha})_+$ is the characteristic function $\widehat{\mu}$ of some symmetric probability distribution $\mu$. Then we have
$$
\widehat{\mu} \Bigl(\frac{t}{n^{1/{\alpha}}} \Bigr)^n = \Bigl( 1 - \frac{|t|^{\alpha}}{n} \Bigr)^n \longrightarrow e^{- |t|^{\alpha}},
$$
which means that the function $e^{- |t|^{\alpha}}$ is also a characteristic function of some $\alpha$-stable probability measure. By the theory of symmetric stable distributions (see e.g. \cite{ST}) we get $\alpha \leqslant 2$. \qed \\

{\bf Example 5.0.} As we have seen in Example 4.0 the probability kernel for the Kingman convolution is equal to the characteristic function $\Phi_s(t) = \Gamma(s+1) \Bigl( \frac{t}{2}\Bigr)^s J_s(t)$ of the variable $\theta_s$ appearing in the definition of this convolution. Consequently, the Kingman convolution $\otimes_{\omega_s}$ is weakly stable for all $s> - \frac{1}{2}$. \\

{\bf Example 5.1.} The classical convolution on $\mathcal{P}_+$ is weakly stable since for its probability kernel $e^{-t}$ we have $g(t) = e^{-|t|}$ which is the characteristic function of the Cauchy distribution. \\

{\bf Example 5.2.} The symmetric convolution is weakly stable since $g(t) = \cos(t)$ is the characteristic function of $\mu_s = \frac{1}{2} \delta_1 + \frac{1}{2} \delta_{-1}$. \\

{\bf Example 5.3.} The $\alpha$-stable convolution $\ast_{\alpha}$ is weakly stable for $\alpha \in (0,2]$ since in this case $e^{-|t|^{\alpha}}$ is the characteristic function of a symmetric $\alpha$-stable measure. For $\alpha>2$ the convolution $\ast_{\alpha}$ is not weakly stable. \\

{\bf Example 5.4.} 
{\bf Example 5.6.} For the Kucharczak convolution the probability kernel is  $\Omega(t) = {{\Gamma(a,t^r)}/{\Gamma(a)}}$ for some $a,r>0$, thus for the function $g(t) = \Omega(|t|)$ we have $g'(t) = - \frac{r}{\Gamma(a)} t^{ar-1} e^{-t^{r}} <0$ for $t>0$ and $g''(t) = \frac{r}{\Gamma(a)} (rt^r +1 -ar) t^{ar-2} e^{-t^r}$, which is positive for $ar\leqslant 1$. This means that the Kucharczak convolution is weakly stable  if $ar\leqslant 1$. \\

{\bf Example 5.7.} For all $n \in \mathbb{N}$ the Kucharczak-Urbanik convolution $\vartriangle_{\alpha, n}$ is weakly stable if $\alpha \in (0,1]$. \\

{\bf Example 5.9.} The Kendall-type convolutions $\vartriangle_{c,\alpha,p}$ with the probability kernel $\varphi_{c,\alpha,p}(t) = (1 - (c+1)t^{\alpha} + ct^{p\alpha}) \mathbf{1}_{[0,1]}(t)$, $p\geqslant 2$, $\alpha>0$, is weakly stable for $\alpha \leqslant 1$ since then, in all admissible cases, $\varphi_{c,\alpha,p}'(t) \leqslant 0$ and $\varphi_{c,\alpha,p}''(t) \geqslant 0$ for all $t \in [0,1]$. This by the P\'olya criterion shows that  $\varphi_{c,\alpha}(|t|)$, $t\in \mathbb{R}$ is a characteristic function of some probability measure $\mu$. This means that $\vartriangle_{c,\alpha,p}$ is $\mu$-weakly stable.

Of course  the $\max$-convolution and $\diamondsuit_{p,\alpha}$ convolution cannot be weak generalized convolutions since they are not regular.

\section{Lack of memory property }
In the classical theory of stochastic processes a very important role plays the Poisson process build on a the sequence of i.i.d.  exponentially distributed random variables. This particular choice of distribution was caused by the lack of memory property exclusively satisfied by the  exponential distribution.  It turns out that a  generalized convolution $\diamond$  admits or not the existence of a distribution with the lack of memory property with respect to $\diamond$. However if such distribution exists, then it is unique up to a scale coefficient. To analyze this notion more precisely we need to define monotonic convolutions first:
\begin{defn}
A generalized convolution $\diamond$ on $\mathcal{P}_+$ is monotonic if
for every $x,y \geqslant 0$ we have
$$
\delta_x \diamond \delta_y \bigl( [x \vee y, \infty) \bigr) =1.
$$
\end{defn}
Informally speaking the generalized convolution is monotonic if the corresponding generalized sum of independent positive random variables cannot be smaller than the biggest of them. \\

{\bf Example 6.0.}
Not every generalized convolution is monotonic. The best known convolution without this property is the Kingman (or Bessel) convolution since for every $s >-\frac{1}{2}$ and $x,y>0$ we have
$$
{\rm supp} \bigl( \delta_x \otimes_{\omega_s} \delta_y \bigr) = \bigl[ |x-y|, x+y \bigr].
$$
\begin{defn}
A probability measure $\nu \in \mathcal{P}_+$ has the lack of memory property with respect to the generalized convolution $\diamond$ if
$$
\mathbf{P} \left\{ X > x\diamond y \big| X >x \right\} = \mathbf{P} \left\{ X>y \right\}, \quad x,y \geqslant 0,
$$
where $X$ is a random variable with distribution $\nu$ and $(x \diamond y)$ is any random variable with  $\mathcal{L}(x \diamond y)  = \delta_x \diamond \delta_y$,  independent of $X$.
\end{defn}
\begin{rem} {\rm
Notice that if the generalized convolution $\diamond$ is monotonic then the equation from Def. 8 can be changed into 
$$
\mathbf{P} \left\{ X > x\diamond y \right\} =  \mathbf{P} \left\{ X>x \right\} \mathbf{P} \left\{ X>y \right\}, \quad x,y \geqslant 0\color{violet}.
$$}
\end{rem}
It was shown in \cite{Poisson}, Prop. 5.2 that the measure $\nu \in \mathcal{P}_+$ with the cumulative distribution function $F$ has the lack of memory property with respect to the monotonic generalized convolution $\diamond$ if and only if the probability kernel $\Omega(t)$ is monotonically decreasing and $F(t) = 1 - \Omega( ct)$, $t>0$,  for some constant $c>0$. In view of the previous considerations we have the following:
\begin{thm}
Let $\diamond$ be a monotonic generalized convolution with the probability kernel $\varphi$. Then the following conditions are equivalent:
\begin{itemize}
\item[1)] $\varphi(t)$ is monotonically decreasing on $\mathbb{R}_+$ and $\varphi(+\infty)=0$;
\item[2)]  $ (1- \varphi(t)) \mathbf{1}_{[0,\infty)}(t)$ is the distribution function of a measure with lack of memory property;
\item[3)] $\varphi(t^{-1}) \mathbf{1}_{[0,\infty)}(t)$ is the cumulative distribution function of some probability measure
\end{itemize}
\end{thm}

\begin{Example}
The classical convolution $\ast$ is evidently monotonic, its probability kernel is $e^{-t}$, thus it admits the distribution with lack of memory property, which is well known to be exponential. 
\end{Example}

\begin{Example} The symmetric convolution $\bowtie$ is not monotonic, since for $x,y>0$
$$
{\rm supp} \bigl(\delta_x \bowtie \delta_y \bigr) = \bigl\{ |x-y|, x+y \bigr\}.
$$
\end{Example}
\vspace{-5mm}
\begin{Example}
The $\alpha$-stable generalized convolution $\ast_{\alpha}$ 
is monotonic and has the kernel of generalized characteristic function $\Omega(t) = e^{-t^{\alpha}}$. This function satisfies assumptions of Thm. 3 thus $\ast_{\alpha}$ admits the distribution with lack of memory property with cumulative distribution function $1- F_Z(t) = e^{-t^{\alpha}} \mathbf{1}_{[0,\infty)}(t)$. The convolution $\ast_{\alpha}$ is $\mu$-weak with respect to $\otimes_{\triangledown}$-convolution, where $\mu$ has the cumulative distribution function $F(t) =1- F_Z(t^{-1})$ and the density
$$
f(t) = \alpha t^{-\alpha -1} e^{-t^{-\alpha}} \mathbf{1}_{(0,\infty)} (t).
$$
\end{Example}
\vspace{-10mm}
\begin{Example}
The Kendall convolution $\vartriangle_{\alpha}$ is monotonic since $\delta_a \vartriangle_{\alpha} \delta_b$, $a,b >0$, is a measure supported in $[ a \vee b, \infty)$ and its probability kernel $\Omega(t) = (1-t^{\alpha})_+$ satisfies the assumptions of Thm. 3, thus the measure $\mu$ with lack of the memory property is ${\rm pow}(\alpha)$ since its cumulative distribution function is $ F(t) = t^{\alpha} \mathbf{1}_{[0,1]}(t) +\mathbf{1}_{[1,\infty)}(t)$ .
\end{Example}

\begin{Example}
The $\max$-convolution is evidently monotonic and its distribution with the lack of memory property is $\delta_1$. Note that the corresponding Poisson process is rather dull as it is not moving at all: $\max\{1,1\}=1 = \max\{1, \max\{1,1\}\}$. 
\end{Example}
\begin{Example}
The Kucharczak convolution for $a\in (0,1]$, $r>0$, is monotonic and its probability kernel is given by $\Omega(t) = \frac{\Gamma(a, t^r)}{\Gamma(a)}$, $t>0$. 
Thus the corresponding distribution with lack of memory property  is the Weibull distribution with the distribution function $F$ such that $F(t) = (1-\Omega(t)) \mathbf{1}_{[0,\infty)}(t)$ and the density
$$
F'(t) = \frac{r}{\Gamma(a)} \, t^{ar-1} e^{-t^{r}} \mathbf{1}_{(0,\infty)} (t).
$$
\end{Example}
\vspace{-8mm}
\begin{Example}
The Kucharczak-Urbanik generalized convolution is monotonic and the function
$$
f(x) = n \alpha t^{\alpha -1} \bigl( 1 - t^{\alpha} \bigr)_+^{n-1}
$$
is the density of its distribution with lack of memory property.
\end{Example}
\begin{Example}
The $\diamondsuit_{p,\alpha}$ generalized convolution is not regular but it is monotonic. It admits the existence of a distribution $\lambda$ with lack of memory property, defined by
$$
\lambda(dx) = \alpha p x^{\alpha -1} dx + (1-p) \delta_1(dx).
$$
\end{Example}
\vspace{-9mm}
\begin{Example}
The Kendall type convolutions are monotonic since their probability kernels $\varphi_{c,\alpha,p}$ are monotonically decreasing. The measure with the lack of memory property has density
$$
f_{c,\alpha,p}(x) = \alpha \bigl[ 1+c - cp x^{\alpha(p-1)} \bigr] x^{\alpha -1} \, \mathbf{1}_{(0,1)}(x).
$$
\end{Example}
 \section{The Kendall convolution expressed by the  $\max$-convolution}
We  can replace the classical convolution in the condition defining weak stability  by any generalized convolution $\diamond$, as it was done by Kucharczak and Urbanik in \cite{KucharczakUrbanik} and by Jasiulis-Go{\l}dyn and Kula in \cite{JasKula}:
\begin{defn}
Let $\diamond$ be a generalized convolution on $\mathcal{P}_+$. A distribution $\mu$ is weakly stable with respect to  $\diamond$ $(\diamond$-weakly stable$)$ if
$$
\forall\, a,b \geqslant 0 \,\, \exists \, \lambda = \lambda_{a,b} \in \mathcal{P}_+ \quad  T_a \mu \diamond T_b \mu = \lambda \circ \mu,
$$
\end{defn}
Distributions  weakly stable with respect to  $\diamond$ define new generalized convolution, called the weak generalized convolution with respect to  $\diamond$.

\begin{defn} 
Let $\mu$ be a weakly stable measure with respect to the generalized convolution $\diamond$. Then a  $\mu$-weak generalized convolution $\otimes_{\mu,\diamond}$ with respect to $\diamond$ is defined as follows: for any $a,b \geqslant 0$ 
$$
\delta_a \otimes_{\mu,\diamond} \delta_b = \lambda \quad \hbox{ if } \quad T_a \mu \diamond T_b \mu = \lambda \circ \mu.
$$
\end{defn}
Equivalently we can say that for every $\lambda_1, \lambda_2, \lambda \in \mathcal{P}_+$
$$
\lambda_1  \otimes_{\mu, \diamond} \lambda_2 = \lambda \quad \hbox{ if } \quad \bigl(\lambda_1 \circ \mu \bigr) \diamond \bigl(\lambda_2 \circ \mu \bigr) = \lambda \circ \mu.
$$
Even though the conditions described in Definitions 9 and 10 suggest a strict  connection between the $\diamond$-weakly stable distribution and $\diamond$-stable distribution   this is not the case. The measure $\lambda$ is $\diamond$ stable if
\begin{equation} \label{wsstrictly}
\forall\, a,b \geqslant 0 \,\, \exists \, c >0, \, \exists\, d \in \mathbb{R} \quad  T_a \lambda \diamond T_b \lambda = T_c \lambda \diamond \delta_{d}.
\end{equation}
If $d = d(a,b) \equiv 0$ then the measure $\lambda$ is called $\diamond$ strictly stable and the generalized characteristic function of $\lambda$ is of the form $ \Phi_{\lambda}^{\diamond} (t) = e^{-A t^{\alpha}}$ for some $A\geqslant 0$ and $\alpha >0$ (see \cite{Urbanik73, Urbanik76, Urbanik88}). 
The $\diamond$-stable measures which are not $\diamond$-strictly stable distributions are studied in a series of papers \cite{JarMis, JJMM, MM, KO}, but we still do not have their complete characterization even in the seemingly easier case of weak generalized convolution. \\

The following Theorem is a continuation of the Thm. 3 describing lack of memory property:
\begin{thm}
Let $\diamond$ be a monotonic generalized convolution with the probability kernel $\varphi$. Then the following conditions are equivalent:
\begin{itemize}
\item[1)] $\diamond$ admits the existence of a distribution with lack of memory property;
\item [4)] $\diamond$ is a weak generalized convolution with respect to the $\triangledown$ convolution based on $\triangledown$-weakly stable measure $\mu$ with the distribution function $\varphi(t^{-1}) \mathbf{1}_{[0,\infty)}(t)$, i.e. $\diamond = \otimes_{\mu, \triangledown}$.
\end{itemize}
\end{thm}

\noindent
{\bf Proof.} Only the implication $1) \rightarrow 4)$ requires explanation: By 3) of Thm. 4 we can consider a random variable $X$ with cumulative distribution function of the form $F_X(t):= \varphi(t^{-1}) \mathbf{1}_{[0,\infty)}(t)$. Since $\varphi \colon [0,\infty) \rightarrow \mathbb{R}$ is the probability kernel of  $\diamond$, then for $a,b >0$
\begin{eqnarray*}
F_{\max\{ aX, bX'\}} (t) & = & F_{aX}(t ) F_{bX}(t) = F_X(ta^{-1} ) F_X(tb^{-1}) \\
 &=& \int_0^{\infty} \varphi(t^{-1}s) \delta_a(ds) \cdot\!\! \int_0^{\infty} \varphi(t^{-1}s) \delta_b(ds) \\
 & = & \int_0^{\infty} \varphi(t^{-1}s) (\delta_a \diamond \delta_b) (ds) = F_{\theta X}(t),
\end{eqnarray*}
where $X'$ is an independent copy of $X$, $\mathcal{L}(\theta) = \delta_a \diamond \delta_b$ and $\theta$ is independent of $X$. \qed

\begin{rem} {\rm
By Thm. 4 we know that the generalized convolution $\diamond$ has a kernel $\Omega$ that is monotonically decreasing to zero iff $\diamond = \otimes_{\mu, \triangledown}$, where $\mu$ is a $\triangledown$-weakly stable probability measure with the cumulative distribution function $F(t) := \Omega(t^{-1}) \mathbf{1}_{[0,\infty)}(t)$ and
\begin{equation}\label{maxstable}
\max\bigl\{ \theta_1 X_1, \theta_2 X_2 \bigr\} \stackrel{d}{=} \theta Z,
\end{equation}
where $\theta, \theta_1, \theta_2$ are i.i.d. with cumulative distribution function $F$, $\mathcal{L}(X_1) \diamond \mathcal{L}(X_2) = \mathcal{L}(Z)$  such that $\theta,\theta_1, \theta_2, X_1, X_2, Z$ are independent. }
\end{rem}
\begin{rem}{\rm
Notice that the measure $\mu$ with cumulative distribution function $F$ is weakly stable with respect to $\triangledown$-convolution if
$$
\forall \, x,y,t > 0 \, \exists \, \lambda \in \mathcal{P}_+ \quad F(xt) F(yt) = \int_0^{\infty} F(st) \lambda(ds).
$$
We do not have here the complete solution of this integral-functional equation but we present a rich list of examples connected with some selected generalized convolutions.
}
\end{rem}
\begin{Example} There is a surprising connection between the classical and $\max$-convolution. The classical convolution $\ast$ on $\mathcal{P}_+$ has the probability kernel $\Omega(t) = e^{-t} \mathbf{1}_{[0,\infty)}(t)$, which satisfies assumptions of Thm. 4. Thus the measure $\mu$ with the cumulative distribution function $F(t) = e^{-t^{-1}} \mathbf{1}_{[0,\infty)}(t)$ and the density
$f(t) = t^{-2} e^{-t^{-1}} \mathbf{1}_{[0,\infty)}(t)$
is $\triangledown$-weakly stable, $\ast = \otimes_{\mu, \triangledown}$ and
\begin{equation}\label{surpr}
\max \bigl\{ \theta_1 X_1,\, \theta_2 X_2 \bigr\} \stackrel{d}{=} \theta_1 \bigl( X_1 + X_2\bigr),
\end{equation}
where $\theta_1, \theta_2$ have distribution $\mu$ and $X_1, X_2$ are arbitrary  non-negative random variables such that $\theta_1, \theta_2, X_1, X_2$ are independent.
\end{Example}
\vspace{-3mm}
\begin{rem} {\rm The equality (\ref{surpr}) is also a simple consequence of the lack of  memory property  of the exponential distribution if we notice that  ${1/\theta_i}$ has the exponential distribution with expectation $1$: For any $u>0$
\begin{align*}
\mathbf{P} \left\{ \theta_1 \bigl(X_1 + X_2\bigr) <u \right\} &=  \mathbf{P} \left\{ \theta_1^{-1} > u^{-1}  \bigl(X_1 + X_2\bigr) \right\} \\
& \stackrel{*}{=}\mathbf{P} \left\{ \theta_1^{-1} > u^{-1} X_1 \right\}\,\mathbf{P} \left\{\theta_2^{-1} > u^{-1}  X_2 \right\} \\ &=\mathbf{P}\left\{\max\{\theta_1X_1,\,\theta_2X_2\}<u\right\},
\end{align*}
where $\stackrel{*}{=}$ follows, upon conditioning with respect to $(X_1,X_2)$, by the lack of memory property of $\theta_1^{-1}$.}
\end{rem}
{\bf Example 7.3.} The stable convolution $\ast_{\alpha}$  has the probability kernel $e^{-t^{\alpha}}$, $\alpha >0$, which satisfies assumptions of Thm. 4. Consequently the measure $\mu$ with the cumulative distribution function  $F(t) = e^{-t^{-1}} \mathbf{1}_{[0,\infty)}$ and density
$$
f(t) = \alpha t^{-\alpha -1} e^{-t^{-\alpha}} \mathbf{1}_{[0,\infty)} (t)
$$
is $\triangledown$-weakly stable and $\ast_{\alpha} = \otimes_{\mu, \triangledown}$. 
This leads to an interesting property: if $\theta_1, \theta_2$ have distributions with the density function $f$, variables $\theta_1, \theta_2, X_1, X_2$ are non-negative and independent then
$$
\max\bigl\{ \theta_1 X_1, \theta_2 X_2 \bigr\} \stackrel{d}{=} \theta_1 \left( X_1^{\alpha} + X_2^{\alpha} \right)^{1/{\alpha}}.
$$
{\bf Example 7.4.} For the Kendall convolution $\vartriangle_{\alpha}$, $\alpha >0$, the probability kernel $(1-t^{\alpha})_+$, $\alpha >0$, satisfies assumptions of Thm. 4 thus $\vartriangle_{\alpha} = \otimes_{\mu, \triangledown}$, where $\mu$ is a measure with the cumulative distribution function $F(t) = ( 1-t^{-\alpha}) \mathbf{1}_{[1,\infty)}(t)$, i.e. $\mu = \pi_{\alpha}$. Consequently:  if $\theta_1, \theta_2$ have distribution $\pi_{\alpha}$, variables $\theta_1, \theta_2, X_1, X_2$ are non-negative and independent then
$$
\max\bigl\{ \theta_1 X_1, \theta_2 X_2 \bigr\} \stackrel{d}{=} \theta_1 \left( X_1 \vartriangle_{\alpha} X_2 \right),
$$
where $\left( X_1 \vartriangle_{\alpha} X_2 \right)$ is any random variable with distribution $\mathcal{L}(X_1) \vartriangle_{\alpha} \mathcal{L}(X_2)$ independent of $\theta_1$. \\

{\bf Example 7.5.} {\rm
Notice that the following, rather trivial, property holds:
$$
\forall \, x,y,t >0 \qquad \mathbf{1}_{[0,1]} (xt) \mathbf{1}_{[0,1]} (yt) = \int_0^{\infty} \mathbf{1}_{[0,1]} (st) \delta_{\max\{x,y\}} (ds).
$$
This means that the cumulative distribution function $F(t) = \mathbf{1}_{[0,1]}(t^{-1})$ corresponds to the measure $\delta_1$, which is weakly stable with respect to the $\max$-convolution. This seems to be interesting, but it is only another way to describe the following, trivial property:
$$
\max \{ \theta_1 X_1, \theta_2 X_2 \} \stackrel{d}{=} \theta_1 \max\{ X_1, X_2\}
$$
for $X_1, X_2, \theta_1, \theta_2$ independent, $\mathcal{L} (\theta_1) = \mathcal{L}(\theta_2) = \delta_1$.} \\

{ \bf Example 7.6.}
The Kucharczak convolution has the probability kernel $\Omega(t) ={{\Gamma(a,t^r)}/{\Gamma(a)}}$ satisfying the assumptions of Thm. 4, thus $F(t) = \Omega(t^{-1}) \mathbf{1}_{[0,\infty)}(t)$ is the cumulative distribution function of a $\triangledown$-weakly stable measure $\mu$ with 
$$
f(t) := F'(t) = \frac{r}{\Gamma(a)} \, t^{-ar-1} e^{-t^{-r}} \mathbf{1}_{(0,\infty)} (t).
$$
Again we have: if $\theta_1, \theta_2$ have distributions with the density function $f$, variables $\theta_1, \theta_2, X_1, X_2$ are non-negative and independent then
$$
\max\bigl\{ \theta_1 X_1, \theta_2 X_2 \bigr\} \stackrel{d}{=} \theta_1 \left( X_1 \circ_1 X_2 \right),
$$
where $\left( X_1 \circ_1 X_2 \right)$ is any random variable with distribution $\mathcal{L}(X_1) \circ_1 \mathcal{L}(X_2)$ independent of $\theta_1$. \\

{\bf Example 7.7.} 
The Kucharczak-Urbanik convolution $\vartriangle_{\alpha,n}$ can be defined by the probability kernel $\Omega_{\alpha,n}(t) = (1-t^{\alpha})_+^n$ and its property: for all $\mu_1, \mu_2 \in \mathcal{P}_+$ there exists $\mu =: \mu_1 \vartriangle_{\alpha,n} \mu_2$ such that 
$$
\int_0^{\infty}\! \Omega_{\alpha,n}(tx) \mu_1(dx) \int_0^{\infty} \!\Omega_{\alpha,n}(ty) \mu_2(dy) = \int_0^{\infty}\! \Omega_{\alpha,n}(tx) \mu(dx).
$$
Evidently the function $\Omega_{\alpha,n}(t)$ satisfies the assumptions of Thm. 4, thus the variable $\theta_n$ with the cumulative distribution function $F_{\alpha, n} (t) = (1-t^{-\alpha})^n \mathbf{1}_{[1,\infty)}(t)$   is weakly stable with respect to the $\max$-convolution $\triangledown$.  Moreover, the Kucharczak-Urbanik  convolution $\vartriangle_{\alpha,n}$ is a weak generalized convolution with respect to $\max$-convolution i.e. $\vartriangle_{\alpha,n} = \otimes_{\mu_n, \triangledown}$ and 
$$
\mathcal{L}(X_1) \vartriangle_{\alpha, n}  \mathcal{L}(X_2) =  \mathcal{L} (Z) \quad \hbox{ iff } \quad  \max\bigl\{ \theta_n X_1, \theta_n' X_2 \bigr\} \stackrel{d}{=}  \theta_n Z,
$$
where $\theta_n, \theta_n'$ are i.i.d. with the distribution $\mu_n$ such that $\theta_n, \theta_n', X_1, X_2, Z$ are independent. It is worth noticing also that if $Q_1, \dots Q_n$ are i.i.d. random variables with Pareto distribution $\pi_{\alpha}$ then
$$
\theta_n \stackrel{d}{=} \max \bigl\{ Q_1, \dots, Q_n \bigr\}.
$$

{\bf Example 7.9.} 
For the Kendall-type generalized convolution $\vartriangle_{c, \alpha, p}$the probability kernel
$$
\varphi_{c,\alpha,p} = \left( 1 - (1+c) t^{\alpha} + ct^{\alpha p} \right) \mathbf{1}_{[0,1]}(t)
$$
is the tail of some cumulative distribution function. 
By Thm. 4 we have that  $\varphi_{c,\alpha,p} (t) \mathbf{1}_{[1,\infty)}(t)$ is the tail of distribution function of a measure with lack of memory property with respect to $\vartriangle_{c,\alpha,p}$ convolution and by Thm. 4 each Kendall type generalized convolution is a $\mu$-weak distribution with respect to the $\max$-convolution $\triangledown$, where $\mu \in \mathcal{P}_+$ has the cumulative distribution function $F(t) := \varphi_{c,\alpha,p} (t^{-1}) \mathbf{1}_{[1,\infty)}(t)$.

\section{Convex linear combination property}

In  this section we give a collection of examples of generalized convolutions with the convex linear combination  property. The generalized Kendall convolution is one of these examples. 

\begin{defn}
The generalized convolution $\diamond$ on $\mathcal{P}_+$ has the convex linear combination property with parameter $n \in \mathbb{N}$, $n \geqslant 2$, if there exist functions $p_0, \dots, p_{n-1} \colon [0,1] \mapsto [0,1]$, $\sum_{k=0}^{n-1} p_k(x) \equiv 1$ and there exist  measures $\lambda_0, \dots, \lambda_{n-1} \in \mathcal{P}_+$ such that
$$
\forall \, x \in [0,1] \quad \quad \delta_x \diamond \delta_1 = \sum_{k=0}^{n-1} p_k(x) \lambda_k.
$$
\end{defn}

{\bf Example 8.4.} 
It is evident that the Kendall convolution has the convex linear combination property with the parameter $n=2$. In fact we know much more, see \cite{MisJas2}: it is the only regular generalized convolution with the convex linear convolution property for $n=2$. \\

{\bf Example 8.5.} The max-convolution (which is not regular) is a trivial example of a generalized convolution with the convex linear combination property with $n=1$, since $\delta_x \triangledown \delta_1 = \delta_1$ for $x\in[0,1]$. \\

{\bf Example 8.7.}\label{KuchU} 
The Kucharczak-Urbanik  convolution  $\vartriangle_{\alpha,n}$, $\alpha>0$, $n\in \mathbb{N}$, is another example of generalized convolution with the convex linear combination property for $n+1$ since by equation (\ref{KUdef})
$$
\delta_x \vartriangle_{\alpha,n} \delta_1 (ds)  
:=(1 - x^{\alpha})^{n} \delta_1(ds) + \sum_{k=1}^n {n \choose k} x^{\alpha k} (1 - x^{\alpha})^{n-k} \mu_{k,n} (ds),
$$
where $\mu_{k,n}$ are probability densities given by (\ref{efkaen}). \\

{\bf Example 8.8.} Every non-regular generalized convolutions $\diamondsuit_{p,\alpha}$, $p\in[0,1]$, $\alpha >0$, described by its probability kernel $\Omega_{\diamondsuit_{p,\alpha}} = ( 1 - p t^{\alpha}) \mathbf{1}_{[0,1]}(t)$ has the convex linear combination property for $n=2$. The $\diamondsuit_{p,\alpha}$-convolution is uniquely determined  for $p \neq \frac{1}{2}$ by 
$$
\delta_x \diamondsuit_{p,\alpha} \delta_1 (ds)= (1-px^{\alpha})\, \delta_1(ds) +px^{\alpha} \,\tfrac{\alpha}{2p-1} \tfrac{2p - s^{q}}{s^{2\alpha+1}} \mathbf{1}_{[1,\infty)}(s) ds,\quad x\in[0,1],
$$
where $q = \frac{\alpha(1-2p)}{(1-p)}$. By continuity, for $p\to {1/2}$ we have
$$
\delta_x \diamondsuit_{1/2,\alpha} \delta_1 (ds)= \bigl(1-\tfrac{x^{\alpha}}{2} \bigr)\, \delta_1(ds) + \tfrac{x^{\alpha}}{2}\, \tfrac{\alpha( 1 + 2 \ln s)}{s^{2\alpha+1}}\mathbf{1}_{[1,\infty)}(s) ds,
$$

{\bf Example 8.9.}
Notice that for Kendall-type generalized convolutions $\vartriangle_{c,\alpha,p}$ in each of the five admissible cases described in \cite{Mej18l-c} we  have $\varphi_{c,\alpha,p}(0) = 1$, $\varphi_{c,\alpha,p}(1) = \varphi_{c,\alpha,p}(+\infty) = 0$ and
$$
\delta_x \vartriangle_{c,\alpha,p} \delta_1 = \varphi_{c,\alpha,p}(x) \, \delta_1 + x^{\alpha p} \, \lambda_1 + (c+1)(x^{\alpha} -x^{\alpha p})\, \lambda_2,
$$
for some probability measures $\lambda_1, \lambda_2 \in \mathcal{P}_+$. This means that it has convex linear combination property with $n=3$.

\section{ Description by random variables}

While constructing stochastic processes with independent increments in the sense of generalized convolution it turns out that we have big trouble if we  study path properties of such processes. This was the reason why the authors of \cite{BJMR} introduced the definition 6.2 of  representability for weak generalized convolutions. Roughly speaking the weak generalized convolution $\diamond$ is representable if there exists a method of unique clear choice of variable $X$ for which $\mathcal{L}(X) =  \mu_1 \diamond \mu_2$. The proper definition of representability of generalized convolution requires more conditions if it is suppose to be used in constructing stochastic processes by their paths - for details see Def. 6.2 in \cite{BJMR}.  \\

For the convenience, we denote by $\theta_1 \diamond \theta_2$ any random variable with distribution $\mathcal{L}(\theta_1) \diamond \mathcal{L}(\theta_2)$ if $\theta_1, \theta_2$  are non-negative and independent . \\

\noindent{\bf Example 9.0.} There are at least three methods of representing the  Kingman convolution $\otimes_{\omega_s}$: \\

\begin{enumerate}
 \item[1)] If $n = 2(s+1) \in \mathbb{N}$ then we are using the weakly stable random vector $\mathbf{U} = (U_1, \dots,U_n)$ with uniform distribution on the unit sphere $S_n$ in $\mathbb{R}^n$. Then for independent random variables $\theta_1, \theta_2$ we choose independent copies  $\mathbf{U}_1, \mathbf{U}_2$  of $\mathbf{U}$ such that $\theta_1, \theta_2, \mathbf{U}_1, \mathbf{U}_2$ are independent. Next we define an adding operator on pairs $(\theta_i, \mathbf{U}_i)$, $i=1,2$ by
\begin{equation}\label{Kingeq1}
\theta_1 \mathbf{U}_1 + \theta_2, \mathbf{U}_2 =  \| \theta_1 \mathbf{U}_1 + \theta_2 \mathbf{U}_2 \|_2 \cdot \frac{\theta_1 \mathbf{U}_1 + \theta_2 \mathbf{U}_2}{\| \theta_1 \mathbf{U}_1 + \theta_2 \mathbf{U}_2\|_2} 
\end{equation}
where $\| \cdot \|_2$ denotes the Euclidean norm in $\mathbb{R}^n$. The two product factors on the right are independent and
$$
 \| \theta_1 \mathbf{U}_1 + \theta_2 \mathbf{U}_2 \|_2 \stackrel{d}{=} \theta_1 \otimes_{\omega_s} \theta_2, \qquad
\frac{\theta_1 \mathbf{U}_1 + \theta_2 \mathbf{U}_2}{\| \theta_1 \mathbf{U}_1 + \theta_2 \mathbf{U}_2\|_2}  \stackrel{d}{=}  \mathbf{U}_1.
$$
We see that the equality (\ref{Kingeq1}) is the equality ($\ast$) given in Section 5, following Definition 5, written in the language of random elements, where $\mu = \mathcal{L}(U)$ and
\begin{eqnarray*}
\theta & = & \theta(\theta_1, \mathbf{U}_1, \theta_2, \mathbf{U}_2) =  \| \theta_1 \mathbf{U}_1 + \theta_2 \mathbf{U}_2 \|_2, \\
\mathbf{U} & = & \mathbf{U}( \theta_1, \mathbf{U}_1, \theta_2, \mathbf{U}_2) = \frac{\theta_1 \mathbf{U}_1 + \theta_2 \mathbf{U}_2}{\| \theta_1 \mathbf{U}_1 + \theta_2 \mathbf{U}_2\|_2}.
\end{eqnarray*} \\

 \item[2)]  Recently Misiewicz and Volkovich showed in \cite{MisiewiczVolodia1} that for arbitrary $s> - \frac{1}{2}$ the random vector $\mathbf{W} = (W_1, W_2)$ with the  density proportional to $(1 - x^2 - y^2)^{s - \frac{1}{2}}$ is weakly stable. Moreover for every choice of independent $\theta_1, \theta_2$, random vectors $\mathbf{W}_1, \mathbf{W}_2$ independent copies of $\mathbf{W}$ such that $\theta_1, \theta_2, \mathbf{W}_1, \mathbf{W}_2$ are independent we have
 \begin{equation}\label{Kingeq2}
\theta_1 \mathbf{W}_1 + \theta_2, \mathbf{W}_2 =  \| \theta_1 \mathbf{W}_1 + \theta_2 \mathbf{W}_2 \|_2 \cdot \frac{\theta_1 \mathbf{W}_1 + \theta_2 \mathbf{W}_2}{\| \theta_1 \mathbf{W}_1 + \theta_2 \mathbf{W}_2\|_2}. 
\end{equation}
The two two product factors on the right are independent and
$$
\| \theta_1 \mathbf{W}_1 + \theta_2 \mathbf{W}_2 \|_2 \stackrel{d}{=} \theta_1 \otimes_{\omega_s} \theta_2, \qquad
\frac{\theta_1 \mathbf{W}_1 + \theta_2 \mathbf{W}_2}{\| \theta_1 \mathbf{W}_1 + \theta_2 \mathbf{W}_2\|_2}  \stackrel{d}{=}  \mathbf{W}_1.
$$
The equality (\ref{Kingeq2}) is the equality from Def. 6 in Section 5, written in the sense of equality almost everywhere and
\begin{eqnarray*}
\theta & = & \theta(\theta_1, \mathbf{W}_1, \theta_2, \mathbf{W}_2) =  \| \theta_1 \mathbf{W}_1 + \theta_2 \mathbf{W}_2 \|_2, \\
\mathbf{W} & = & \mathbf{W}( \theta_1, \mathbf{W}_1, \theta_2, \mathbf{W}_2) = \frac{\theta_1 \mathbf{W}_1 + \theta_2 \mathbf{W}_2}{\| \theta_1 \mathbf{W}_1 + \theta_2 \mathbf{w}_2\|_2}.
\end{eqnarray*}
Notice that $\mathbf{W}$ can be identified with the vector $(\cos \phi, \sin \phi)$, where $\phi$ is a random variable with the density proportional to $(\sin^2\varphi)^{s+ \frac{1}{2}}$ on the interval $[0,2\pi]$. Moreover, the vector  $\mathbf{W}$  is living on the unit sphere in $\mathbb{R}^2$, but it does not have uniform distribution there. \\

\item[3)] For any $s> - \frac{1}{2}$ Kingman  in \cite{King} gave the following explicit formula for the random variable $\theta_1 \otimes_{s} \theta_2$:
$$
\theta_1 \otimes_{s} \theta_2 \stackrel{d}{=} \sqrt{\theta_1^2 + \theta_2^2 + 2 \theta_1 \theta_2 \cos\phi},
$$
where $\phi$ is a random variable with the density proportional to the function $(\sin^2\varphi)^{s+ \frac{1}{2}}$ on the interval $[0,2\pi]$. It is known (and easy to check) that if $\phi_1, \phi_2$ are independent copies of $\phi$ then $\cos(\phi_1 - \phi_2) \stackrel{d}{=} \cos \phi$.

This leads to the following Kingmam's interpretation: if $Q$ is a vector of the length $\theta$ forming the angle $\varphi$ with the fixed straight line then we will use the notation $Q = (\theta, \cos \varphi)$. Consequently, using elementary geometry  we have 
$$
\bigl(\theta_1, \cos \varphi_1\bigr) \oplus \bigl(\theta_2, \cos \varphi _2 \bigr)   \stackrel{def}{=}  \Bigl( \sqrt{\theta_1^2 + \theta_2^2 + 2 \theta_1 \theta_2 \cos(\varphi_1 - \varphi_2)}, \cos(\varphi_1 - \varphi_2)\Bigr),
$$
and, by the previous considerations, 
$$
\sqrt{\theta_1^2 + \theta_2^2 + 2 \theta_1 \theta_2 \cos(\phi_1 - \phi_2)} \stackrel{d}{=} \theta_1 \otimes_{s} \theta_2.
$$
\end{enumerate}

In view of Example 9.0. we see that the random variable  $\theta_1 \diamond \theta_2$ with the distribution $\mathcal{L}(\theta_1) \diamond \mathcal{L}(\theta_2)$ can be expressed in many different ways. If we want to base on this representation the construction of stochastic processes with independent (with respect to the generalized convolution $\diamond$) increments only the first representation Example 9.0.1 is admissible - for details see \cite{BJMR}. 

\begin{thm} 
If there exists a function $\psi \colon \mathbb{R}^2 \mapsto \mathbb{R}$ such that
$$
\Psi(\theta_1, \theta_2) (\omega) = \psi(\theta_1(\omega), \theta_2(\omega)) \quad a.e.
$$
for all independent $\theta_1, \theta_2$ then there exists $\alpha \in (0, \infty]$ such that 
$$ 
\psi(x,y) = \bigl( |x|^{\alpha} + |y|^{\alpha} \bigr)^{1/{\alpha}}, \qquad x,y \in \mathbb{R}^2,
$$
which follows from the Bohnenblust theorem (for details  see \cite{BJMR}).
\end{thm}

\noindent
Almost trivially we have the following representations of discussed here convolutions by random variables:
$$
\mathcal{L}({\theta_1}) \ast \mathcal{L}({\theta_2})  = \mathcal{L}(\theta_1 + \theta_2). \leqno{{\bf Example\,\, 9.1.}}
$$
$$
\mathcal{L}({\theta_1}) \bowtie \mathcal{L}({\theta_2})  = \mathcal{L}\bigl((\theta_1 + \theta_2) Q + |\theta_1 - \theta_2| (1-Q)\bigr), \leqno{{\bf Example\,\, 9.2.}}
$$
where $\mathbf{P}\{ Q=1\} = \mathbf{P}\{ Q = 0\} = \frac{1}{2}$ such that $Q, \theta_1, \theta_2$ are independent. 
$$
\mathcal{L}({\theta_1}) \ast_{\alpha} \mathcal{L}({\theta_2})  = \mathcal{L}( (\theta_1^{\alpha} + \theta_2^{\alpha})^{1/{\alpha}}). \leqno{{\bf Example\,\, 9.3.}}
$$
$$
\mathcal{L}({\theta_1}) \triangledown  \mathcal{L}({\theta_2})  = \mathcal{L}(\max\{\theta_1, \theta_2\}). \leqno{{\bf Example\,\, 9.5.}}
$$
\begin{thm}
Assume that the generalized convolution $\diamond$ on $\mathcal{P}_+$  has the convex linear combination property. Then $\diamond$ is represented by random variables.
\end{thm}

{\bf Proof.} Assume that $\mathcal{L}(\theta_1) = \mu_1$, $\mathcal{L}(\theta_2) = \mu_2$ such that $\theta_1, \theta_2$ are independent.  By our assumptions for every $x\in [0,1)$ there exist $n\in \mathbb{N}$, $p_0, \dots, p_{n-1} \colon [0,1] \mapsto [0,1]$, $\sum_{k=0}^{n-1} p_k(x) = 1$ for all $x\in [0,1]$, and there exist measures $\lambda_0, \dots, \lambda_{n-1} \in \mathcal{P}_+$ such that
\vspace{-2mm}
\begin{equation} \label{lcp}
\forall\, x \in [0,1] \qquad \delta_x \diamond \delta_1 = \sum_{k=1}^{n-1} p_k(x) \, \lambda_k.
\end{equation}
Now we define some auxiliary random variables: $M= \max \{\theta_1, \theta_2\}$, $m = \min \{ \theta_1, \theta_2 \}$ and $\varrho = \varrho(\theta_1, \theta_2) :=  {m/M}$. For the numbers $s_0(x) = p_0(x)$, $s_k(x) = \sum_{j=0}^{k-1} p_j(x)$ for $k=1, \dots, n-1$, we define a sequence of intervals:   $A_0(x) = [0, p_0(x)]$ and
$$
A_k(x) = \bigl(s_{k-1}(x), s_{k}(x)\bigr], \quad k=1,\dots, n-1.
$$
Of course $\bigcup_{k=0}^{n-1} A_k(x) = [0,1]$ for all $x \in [0,1]$. Now we choose random variables $Q_0, \dots, Q_{n-1}$ with distributions $\lambda_0, \dots, \lambda_{n-1}$ respectively, a random variable $U$ with uniform distribution on  the interval $[0,1]$ such that $\theta_1, \theta_2, \theta_3, U, Q_0, \dots, Q_{n-1}$ are independent. Now we are able to define the random variables representing the convolution $\lambda_1 \diamond \lambda_2$:
$$
x \diamond 1 \stackrel{d}{=} \sum_{k=0}^{n-1} \mathbf{1}_{A_k(x)} (U) Q_k,
$$
and
$$
\lambda_1 \diamond \lambda_2 = \mathcal{L} \Biggl( M \sum_{k=0}^{n-1} \mathbf{1}_{A_k(\varrho)} (U) \,Q_k \Biggr). \eqno{\Box}
$$

\noindent
{\bf Example 9.4.} For representability of the Kendall convolution take non-negative independent random variables $\theta_1, \theta_2$ and  we define, as in the proof of Thm. 5, $M= \max \{\theta_1, \theta_2\}$, $ m = \min\{ \theta_1, \theta_2\}$,  $\varrho ={m/M}$. 
Let $U$ has the uniform distribution on $[0,1]$, $\Pi_{2\alpha}$ has the Pareto distribution $\pi_{2\alpha}$ and $U$,  $\Pi_{2\alpha}$ and $\theta_1,\theta_2$ are independent. Then
$$
\theta_1\vartriangle_{\alpha} \theta_2 \stackrel{d}{=}  M \bigl(\mathbf{1}_{(\varrho^{\alpha},1]}(U) + \Pi_{2\alpha} \mathbf{1}_{[0,\varrho^{\alpha}]} (U) \bigr). 
$$
Another representation of $\theta_1\vartriangle_{\alpha} \theta_2$, found in \cite{bigone} or directly obtained from Thm. 1. Since $\mathbf{P} \{ \frac{\theta_i}{Z_i} < t\} = G_i(t)$, we have the following:
$$
\theta_1\vartriangle_{\alpha} \theta_2 \stackrel{d}{=} \max \left\{ \max\{ \theta_1, \theta_2\}, \min\left\{ \frac{\theta_1}{Z_1}, \frac{\theta_2}{Z_2} \right\} \right\},
$$
where $Z_1, Z_2$ are i.i.d. with ${\rm pow}(\alpha)$ distribution such that $\theta_1, \theta_2, Z_1, Z_2$ are independent.
\vspace{-2mm}
\begin{rem} {\rm
The construction proposed in Thm. 6 can be trivially adapted to Examples 9.7, 9.8 and 9.9, thus we have that the Kucharczak-Urbanik convolutions, $\diamondsuit_{p,\alpha}$-convolutions and Kendall type convolutions can be represented by random variables. }
\end{rem}

\noindent
{\bf Example 9.7.} For the Kucharczak-Urbanik convolution representation by random variables can be done in a more interesting way: 

We introduce first an useful notation: for any $1\leqslant k\leqslant n$ define a function $\sigma_{k,n}:\mathbb R^n\to\{1,\ldots,n\}$ by
$$\sigma_{k,n}(x_1,\ldots,x_n)=x_j\quad \Leftrightarrow\quad \#\{i\in\{1,\ldots,n\}:\,x_i\leqslant x_j\}=k,
$$
for any $j, k \in\{1,\ldots,n\}$. If $X_1,\ldots,X_n$ are i.i.d. random variables the random variable $X_{k:n} := \sigma_{k,n}(X_1,\ldots,X_n)$ is called the $k$'th order statistics (based on $n$ i.i.d. observations), $k=1,\ldots,n$. In particular, $X_{1:n}=\min\{X_1,\ldots,X_n\}$ and $X_{n:n}=\max\{X_1,\ldots,X_{n:n}\}$.  For basic information on order statistics see e.g. \cite{DavNag,WA2004}.  \\

We need also to notice that if $Q$ is the Pareto random variable with distribution $\pi_{\alpha}$, then $Q^{-1}$ has the power distribution $\mathrm{pow}(\alpha)$ with the density $\alpha x^{\alpha-1} \mathbf{1}_{[0,1]}(x)$. Moreover,  if $V_i=Q_i^{-1}$, $i=1,\ldots,n$, are i.i.d. variables with the  power distribution $\mathrm{pow}(\alpha)$ then  
$$
Q_{k:n}=V_{n-k+1:n}^{-1}\quad k=1,\ldots,n.
$$

\begin{thm}\label{KUrep}
Let $\theta_1$ and $\theta_2$ be independent non-negative random variables with distributions $\mu_1$ and $\mu_2$. Then  $\mu_1\vartriangle_{\alpha,n}\mu_2$ is the distribution of the random variable 
$$
M(\theta_1, \theta_2)\,\sum_{k=0}^n\,Q_{k:n+k}\,{\bf 1}_{\bigl(W_{k:n},\,W_{k+1:n}\bigr]}\bigl(\varrho(\theta_1,\theta_2)\bigr),
$$
where $Q_1,\ldots,Q_{2n}$ are i.i.d. random variables with the Pareto distribution $\pi_{\alpha}$, $W_1, \dots, W_n$ are i.i.d. random variables with the distribution ${\rm pow}(\alpha)$ such that $Q_1,\ldots,Q_{2n}, W_1, \dots, W_n$ are independent and $Q_{0:n}:=1, W_{n+1:n}=\infty$.
\end{thm}
{\bf Proof.}
Note that the basic components of the Kucharczak-Urbanik convolution, see \eqref{KUdef}, are probability measures with the densities $f_{k,n}$, $n \in \mathbb{N}$, $k =1,\dots,n$,  defined in \eqref{efkaen}. The key observation here is that $f_{k,n}$ is the density of $Q_{k:n+k}$ where $Q_1,\ldots,Q_{2n}$ is an i.i.d. sample from the same Pareto $\pi_{\alpha}$   distribution. Now by \eqref{KUdef} in Section 3 we have:
$$
x \vartriangle_{\alpha,n} 1 \stackrel{d}{=} \sum_{k=0}^n Q_{k:n+k} {\large \mathbf{1}}_{\{ B_n(x^{\alpha}) = k\}},
$$
where $B_n(x^{\alpha})$ is the Bernoulli random variable (counting successes in $n$ trials with the success probability $p=x^{\alpha}$) such that $B_n(x^{\alpha})$ and $(Q_1,\ldots,Q_{2n})$ are independent. 

It remains to show that for all $k=0,1,\dots, n$ we have
$$
\mathbf{P} \bigl\{ B_n(x^{\alpha}) = k \bigr\} = \mathbf{E}\, \mathbf{1}_{(W_{k:n},W_{k+1:n}]} (x) = \mathbf{P} \bigl\{ W_{k:n}< x \leqslant W_{k+1:n} \bigr\},
$$
where $W_1, \dots, W_n$ are i.i.d. random variables with the distribution ${\rm pow}(\alpha)$. To see this we recall (see e.g. \cite{DavNag}) that the bivariate density function  $f_{k,k+1:n}$ of $(X_{k:n},X_{k+1:n})$ for i.i.d. random variables  $X_1,\ldots,X_n$ with the density $f$ and cumulative distribution function  $F$ has the form
$$
f_{k,k+1:n}(x,y)=\frac{n!}{(k-1)!(n-k-1)!}F^{k-1}(x)F^{n-k-1}(y)f(x)f(y){\large\bf 1}_{\{x<y\}}.
$$
Therefore, for any $r$
\begin{align*}
& P\bigl\{X_{k:n}<r\leqslant X_{k+1:n} \bigr\} = \\ & \,\frac{n!}{(k-1)!(n-k-1)!} \int_{-\infty}^r\!\! F^{k-1}(x) f(x)\,dx \int_r^{\infty} \! (1-F(y))^{n-k-1} f(y)\,dy\\ & = \,  \binom{n}{k}F^k(r)(1-F(r))^{n-k}=\mathbf{P} \bigl\{ B_n(F(r)) = k \bigr\}.
\end{align*}
The last formula applied to $W_{k:n},\,W_{k+1:n}$ yields $P\{W_{k:n}<r\leqslant W_{k+1} \} = \mathbf{P} \bigl\{ B_n(x^{\alpha}) = k \bigr\}$. Now, assuming that $Q_1, \dots, Q_{2n}$ and $W_1,\dots, W_n$ are independent, we have
\begin{equation}\label{repre1}
x \vartriangle_{\alpha,n} 1 \stackrel{d}{=}  \sum_{k=0}^n\,Y_{k:n+k} \mathbf{1}_{\bigl(W_{k:n},W_{k+1:n}\bigr]} (x). 
\end{equation}
In order to get the final statement it is enough to choose $Q_1, \dots, Q_{2n}$ and $W_1,\dots, W_n$ independent of $\theta_1, \theta_2$ and notice that
$$
\theta_1 \vartriangle_{\alpha,n} \theta_2 =  M(\theta_1, \theta_2) \Bigl( \delta_{\varrho(\theta_1, \theta_2)} \vartriangle_{\alpha,n} \delta_1 \Bigr).  \eqno{\Box}
$$
\vspace{-5mm}
\begin{rem} {\rm
Notice that for the generalized convolution $\diamond$ on $\mathcal{P}_+$  with the convex linear combination property we have
\begin{eqnarray*}
\frac{1}{\theta_1 \diamond \theta_2} & \stackrel{d}{=} & \frac{1}{M(\theta_1, \theta_2) \sum_{k=0}^{n-1} \mathbf{1}_{A_k(\varrho(\theta_1, \theta_2))} (U) \, X_k} \\
& = & m(\theta_1^{-1}, \theta_2^{-1}) \sum_{k=0}^{n-1} \mathbf{1}_{A_k(\varrho(\theta_1, \theta_2))} (U) \, X_k^{-1},
\end{eqnarray*}
if $\theta_1, \theta_2, X_0, \dots X_{n-1}$ are independent, $\mathcal{L}(X_k) = \lambda_k$, $k=0,\dots, n-1$ as in the representation \eqref{lcp}. We used here equality $\varrho(\theta_1, \theta_2) = \varrho(\theta_1^{-1}, \theta_2^{-1})$.}
\end{rem}

\begin{rem} {\rm
Applying this techniques to the Kucharczak-Urbanik convolution   $\vartriangle_{\alpha,n}$ and using the result of Theorem 1 we obtain
$$
\frac{1}{\theta_1 \vartriangle_{\alpha,n} \theta_2}\, \stackrel{d}{=}\, m(\theta_1, \theta_2)\,\sum_{k=0}^n\,{\bf 1}_{\bigl(W_{k:n},\,W_{k+1:n}\bigr]}\bigl(\varrho(\theta_1,\theta_2)\bigr) \,V_{n+1:n+k}^{-1},
$$
where  $V_1, \dots, V_n, W_1, \dots, W_n$ are i.i.d. random variables with the distribution ${\rm pow}(\alpha)$ such that $V_{0:n}:=1$, $W_{n+1:n}=\infty$.}
\end{rem}

{\bf Acknowledgements.} This paper is partially supported by the  project "First order Kendall maximal autoregressive processes and their applications", Grant no POIR.04.04.00-00-1D5E/16, which is carried out within the POWROTY/ REINTEGRATION programme of the Foundation for Polish Science co-financed by the European Union under the European Regional Development Fund. Jacek Weso{\l}owski was partially supported by the grant 2016/21/B/ST1/00005 of National Science Centre, Poland.

\normalsize

\end{document}